\documentclass[12pt,draftcls,onecolumn]{IEEEtran}
\usepackage{graphics} % for pdf, bitmapped graphics files
\usepackage{epsfig} % for postscript graphics files
\usepackage{amsmath}
\usepackage{array}
\usepackage{mdwtab}
\usepackage{cite}
\usepackage{multirow}
\usepackage{threeparttable}
\usepackage{amsmath,amssymb}
\usepackage{latexsym}
\usepackage{graphicx}
\usepackage{epstopdf}
\usepackage{times} % assumes new font selection scheme installed
\usepackage{balance}
\usepackage{breqn}
\usepackage{threeparttable}
\usepackage{mathtools, cuted}
\usepackage{lipsum, color}
\usepackage{multicol}
\usepackage{tikz}
\usepackage{standalone}
\usepackage{tabularx}
\usetikzlibrary{shapes,arrows}
\usetikzlibrary{automata,arrows,positioning,calc}

 \DeclareMathAlphabet{\mathcal}{OMS}{cmsy}{m}{n}

\newtheorem{remark}{{\bf Remark}}
\newtheorem{lemma}{{\bf Lemma}}
\newtheorem{theorem}{{\bf Theorem}}
\newtheorem{asp}{{\bf Assumption}}
\newtheorem{corol}{{\bf Corollary}}

\usepackage{fancyhdr}
\usepackage[nodayofweek,USenglish]{datetime} 
\IEEEoverridecommandlockouts

\begin{document}
\usdate
\title{Stability analysis of event-triggered anytime control with multiple control laws}
\author{Thuy~V.~Dang,~\IEEEmembership{}
        K.~V.~Ling,~\IEEEmembership{Member,~IEEE,}
        and~Daniel~E.~Quevedo,~\IEEEmembership{Senior~Member,~IEEE}% <-this % stops a space
%\thanks{Manuscript received May 22, 2017; first revision December 19, 2017; second revision March 5, 201.}
\thanks{Thuy V. Dang is with the Interdisciplinary Graduate School (IGS) and with the School of Electrical and Electronic Engineering, Nanyang Technological University, 639798 Singapore. {DANG0028@e.ntu.edu.sg}}% <-this % stops a space
\thanks{K.V. Ling is with the School of Electrical and Electronic Engineering, Nanyang Technological University, 639798 Singapore. {EKVLING@ntu.edu.sg}}% <-this % stops a space
\thanks{D.E. Quevedo is with the Faculty of Electrical Engineering and Information Technology (EIM-E), Paderborn University, Warburger Str. 100, 33098 Paderborn, Germany. {dquevedo@ieee.org}}
\thanks{Thuy V. Dang and K.~V. Ling acknowledge the support of the National Research Foundation, Prime Ministers Office, Singapore under its CREATE programme.}
\thanks{\textcopyright 20xx IEEE. Personal use of this material is permitted. Permission from IEEE must be obtained for all other uses, in any current or future media, including reprinting/republishing this material for advertising or promotional purposes, creating new collective works, for resale or redistribution to servers or lists, or reuse of any copyrighted component of this work in other works.
}
}
%\markboth{Submitted to IEEE Transactions on {Automatic Control}}%
\maketitle
\tikzstyle{blockrec} = [draw, fill=blue!20, rectangle, 
    minimum height=3em, minimum width=6em]
\tikzstyle{ctrlBlk} = [draw, fill=blue!20, rectangle, 
    minimum height=0.4cm, minimum width=6em]
\tikzstyle{blocksquare} = [draw, fill=blue!20, rectangle, 
        minimum height=4em, minimum width=4em]
\tikzstyle{blocklong} = [draw, fill=blue!20, rectangle, 
                minimum height=1em, minimum width=6em]
\tikzstyle{eventblk} = [draw, fill=gray!20, rectangle, 
                                minimum height=0.5cm, minimum width=6cm]
\tikzstyle{smallCir} = [draw, fill=black!20, circle, node distance=1cm]
\tikzstyle{smallsquare} = [draw, fill=black, rectangle,minimum height=0.1cm,minimum width=0.1cm]
\tikzstyle{input} = [coordinate]
\tikzstyle{output} = [coordinate]
\tikzstyle{pinstyle} = [pin edge={to-,thin,black}]
\setlength{\abovedisplayskip}{3pt}
\setlength{\belowdisplayskip}{3pt}

\begin{abstract}
To deal with time-varying processor availability and lossy communication channels in embedded and networked control systems, one can employ an event-triggered sequence-based anytime control (E-SAC) algorithm. The main idea of E-SAC is, when computing resources and measurements are available, to compute a sequence of tentative control inputs and store them in a buffer for potential future use. State-dependent Random-time Drift (SRD) approach is often used to analyse and establish stability properties of such E-SAC algorithms. However, using SRD, the analysis quickly becomes combinatoric and hence difficult to extend to more sophisticated E-SAC. In this technical note, we develop a general model and a new stability analysis for E-SAC based on Markov jump systems. Using the new stability analysis, stochastic stability conditions of existing E-SAC are also recovered. In addition, the proposed technique systematically extends to a more sophisticated E-SAC scheme for which, until now, no analytical expression had been obtained. 
\end{abstract}
%\IEEEpeerreviewmaketitle
\section{Introduction}
It is common in the embedded or networked control system that processor availability varies due to varying computational loads and multi-tasking operations. Anytime control algorithm was first proposed in \cite{bhattacharya2004anytime} to deal with time-varying processor availability.
It uses an idea from the AI community called {\it anytime algorithm}\cite{zilberstein1996using}, which is a computational procedure that could provide a valid answer even when it is terminated prematurely.

There are various forms of anytime control algorithms. In \cite{bhattacharya2004anytime}, different controllers with different floating point operations were designed. Notable later works include \cite{greco2011design,gupta2013control,quevedo2013sequence,pant2015codesing}. In \cite{greco2011design}, a stochastic switching law within a set of pre-designed controllers is proposed. In \cite{gupta2013control}, the main idea is to sequentially calculate the components of the plant input vector. In \cite{quevedo2013sequence}, named as sequence-based anytime control (SAC), a buffer is used to store the tentative future inputs. In \cite{pant2015codesing}, a method for co-design of estimator and controller is proposed where the controller requests a time-varying criterion for the estimator.  When sensor measurements are transmitted through a communication network, the measurements may be unavailable due to packet dropouts, or network congestion. Among anytime control algorithms, SAC can handle this situation since it has a buffer which serves to provide a control input even when no measurement is received.

Motivated by the idea of using event-triggered control (see e.g. \cite{aastrom2002comparison,tabuada2007event,borgers2017riccati,wang2017dynamic,ma2017event}) as a method to reduce demands on the network and computing processor while guaranteeing satisfactory levels of performance\cite{postoyan2015framework}, SAC with an event-triggering mechanism (E-SAC) was proposed in \cite{quevedo2014stochastic} and the  State-dependent Random-time Drift (SRD) technique of \cite{yuksel2013random} was employed to analyse the stability of E-SAC.
 
In our conference contribution \cite{huang2014event}, the E-SAC was extended to a more sophisticated scheme featuring two control laws, a coarse and a fine law. The fine control law could be viewed as an improved version of the coarse control law that requires more processing resources than the coarse control law. Such ideas are wide-spread, e.g., in Model Predictive Control (MPC) \cite{maciejowski2002predictive,ling2012multiplexed,mayne2014model,grune2017nonlinear}, to compute sub-optimal and optimal solutions are two strategies that one can choose depending on available computation time. Alternatively, fixed-point and floating-point implementations can be used for trading off computation time and quality (accuracy) \cite{constantinides2011numerical}. 
 
It was demonstrated in \cite{huang2014event} that with the multi-control law E-SAC schemes, the communication and processing resources could be used more efficiently. Performance in terms of empirical closed-loop cost, channel utilisation and regions for stochastic stability guarantees could be improved, compared with the basic E-SAC.

In \cite{huang2014event}, the SRD technique was used to analyse the stability of the proposed multi-control law E-SAC schemes. Unfortunately, this requires one to list all possibilities and the corresponding probabilities. For example, in the two-control law schemes, there are two random variables: (1) the number of times each control law is active during (2) the time interval that the buffer becomes empty again. Therefore, it is a combinatoric problem and quickly becomes intractable.
As a result, a closed-form expression for stability condition cannot be readily obtained by the SRD approach. It was concluded that extending SRD technique to more sophisticated E-SAC schemes will be difficult.

In the present work, we propose a new approach to investigate the stochastic stability of E-SAC schemes. By modelling E-SAC as a Markov jump system (MJS) \cite{ji1990jump} with event-triggering, assuming that processor availability and packet dropouts are identical independent distributed (i.i.d) random processes, we systematically establish stochastic stability guarantee of both one- and multi-control law E-SAC schemes. Our proposed approach recovers stability conditions of \cite{quevedo2014stochastic} which is a one-control law E-SAC scheme.

The remainder of this paper is organised as follows: In Section II we provide a review of the basic E-SAC scheme, including the one control law as proposed in \cite{quevedo2014stochastic}, and the multi-control law E-SAC schemes as proposed in \cite{huang2014event}. In Section III, we propose the Markov jump system with event-triggering (E-MJS) model for stability analysis of the multi-control law schemes. Section IV investigates stochastic stability issues of this E-MJS model. Section V presents the stochastic stability results of E-SAC, derived by the new approach. Section VI documents a simulation study. Section VII draws conclusions.

\textit{Notation:}
$\mathbb{N}=\{1,2,...\}$ represents natural numbers, $\mathbb{N}_0\triangleq \mathbb{N} \cup \{0\}$; $\mathbb{R}$ represents real numbers, $\mathbb{R}_{\geq 0} \triangleq [0,+ \infty)$, $\mathbb{R}_{> 0} \triangleq (0,+ \infty)$; $\{x\}_{\mathcal{K}}$ stands for $\{x_k:~k\in \mathcal{K} \},~\mathcal{K} \subseteq \mathbb{N}_0$. $\sigma(M)$ denotes the spectral radius of matrix $M$. $|x|=\sqrt{x^Tx}$ denotes the Euclidean norm of vector $x$. A function $\phi:~\mathbb{R}_{\geq 0} \rightarrow \mathbb{R}_{\geq 0}$ is of $class-\mathcal{K}_{\infty}$ ($\phi \in \mathcal{K}_{\infty}$) if it is continuous, zero at the origin, strictly increasing and unbounded. $\textbf{Pr}\{\Omega\},~\textbf{Pr}\{\Omega|\Gamma\}$ denote the probability of an event $\Omega$, and the conditional probability of $\Omega$ given $\Gamma$ respectively. The expected value of a random variable $\nu$ given $\Gamma$ is denoted by $\textbf{E}\{\nu|\Gamma\}$, and $\textbf{E}\{\nu\}$ represents the unconditional expectation. For a vector $y$, $y \succ 0$ means that all of its elements are positive. For a matrix $A\in \mathbb{R}^{n\times n}$, $A_{[i:j;l:m]}$ denotes a block matrix contained in $A$ whose elements are taken from row $i$ to $j$, and column $l$ to $m$ of $A$. For $z\in \mathbb{R}$, $\lfloor z \rfloor$ denotes the largest integer that is not bigger than $z$. For $n,k\in \mathbb{N}$, $n ~\text{(mod)}~ k$ means the remainder of $n$ divided by $k$; ${\bf 0}_n$ is a zero vector with dimension $n$. For a vector $x$, $x(i)$ ($i\in \mathbb{N}$) denotes the i-th element of $x$. For a matrix $A$, $||A||_{\infty}$ denotes the infinity norm of $A$.

\section{Review: Event-triggered sequence-based anytime control (E-SAC) schemes}

\begin{figure}[t!]
\begin{tikzpicture}[auto, node distance=2cm,>=latex']
    % We start by placing the blocks
   % \node [input, name=input](input) {(processing power)} ;
    \node (input) {Processor availability} ;
    \node [blocksquare, below of=input, node distance=1.5cm, align=center] (processor) {E-SAC};
    \node [blockrec, right of=processor,pin={[pinstyle]above:$w_k$},node distance=3cm,align=center] (plant) {Plant\\ $x_k$};
    \node [smallCir,right of=plant,node distance=2cm] (s1) {};
    \node [ above right of=s1, node distance=1.25cm] (s2) {$\vert x_k\vert\leq d?$};      
    \node [smallCir,right of=s1,
            node distance=0.75cm,align=center] (event) {};
    \node[output][above of=s1,node distance=0.3cm] (sw3) {};
    \node[output][left of=event,node distance=0.25cm] (sw4) {};
    \draw[<->,bend right] (sw4) to (sw3);
    %\node[above of=s2, node distance=1cm] (eventLabel) {Event $|x_k<d?|$};
    % We draw an edge between the controller and system block to 
    % calculate the coordinate u. We need it to place the measurement block
    \draw [-] (plant) -- (s1);
   \draw [->,line width=1.5pt] (s1) --node{} (s2); 
   % \path (s1)[draw=black,solid,line width=2mm,fill=black,
    %preaction={-triangle 90,thin,draw,shorten >=-1mm}
    %] (s2);
  %  \draw [->] (s2) -- node[name=u] {$x_k$} (event);
    \node [output, right of=event,node distance=1cm] (output) {};
    \node [blocklong, below of=s1] (measurements) {Erasure Channel};
    \node [coordinate, above of=measurements,node distance=0.5cm] (fooChan) {};
    \node[left of=fooChan,node distance=0.45cm] {{\small ${\bf Pr}\{\gamma_k=1\mid\vert x_k\vert> d\}=q$}};
    % Once the nodes are placed, connecting them is easy. 
    \draw [draw,->] (input) -- node {$N_k$ } (processor);
    \draw [->] (processor) -- node {$u_k$} (plant);
    \draw[-] (event)--node{} (output);
    \draw[->] (output)|- (measurements);
   % \draw [->] (system) -- node [name=y] {$y$}(output);
 %   \draw [->] (y) |- (measurements);
    \draw [->] (measurements) -| node[pos=0.99] {} 
        node [near end,right] {$(x_k,\gamma_k)$} (processor);
\end{tikzpicture}
\caption{System Model}
\label{systemModel}
\end{figure}
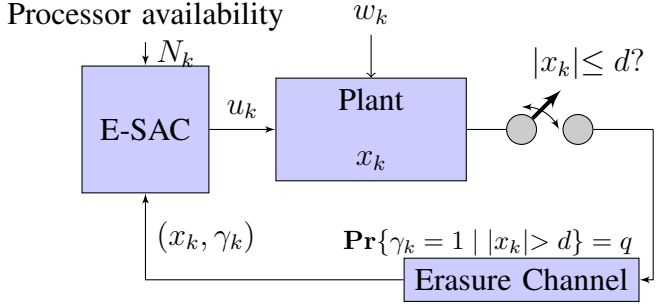

We consider an input-constrained discrete-time non-linear plant model with dynamics given by: 
\begin{equation}
\label{sysmodel}
x_{k+1} = f(x_k,u_k)  \ \ \ \   
\end{equation} where $x_k \in \mathbb{R}^n, ~u_k \in \mathbb{U} \subseteq \mathbb{R} ,~  k \in \mathbb{N}_0$, see Fig. \ref{systemModel}. 

Sensor measurements are transmitted to the controller via a delay-free communication link which introduces packet dropouts. The transmission is a threshold-based event-triggering, i.e., the sensor transmits the measurement only when $|x_k|>d$ where the threshold $d$ is a design parameter. The threshold $d$ is fixed once the system runs, and the triggering event is checked periodically at every sampling instant.

The outcome of the transmission is indicated by the random process ${\{ \gamma \}}_{\mathbb{N}_0}$:  
\begin{equation}
\gamma_k = \left \{
\begin{aligned}[l]
0 & \ \ \mbox{if}\ \ \mbox{the sensor transmits but a dropout occurs} \footnotemark[1]   \\
1 & \ \ \mbox{if}\ \ x_k~\mbox{is received succesfully}   \\
2 & \ \ \mbox{if the sensor did not transmit}\ \ (\mbox{i.e.}~ |x_k| \leq d) 
\end{aligned} \right. \notag
\end{equation}
which is assumed to be (conditionally) independent and identical (i.i.d) with a successful transmission probability
\footnotetext[1]{We assume that packet dropout ($\gamma_k=0$) is distinguishable from no transmission ($\gamma_k=2$), e.g., through error-detection coding and monitoring of received energy/waveforms in the sensor transmission band (see \cite{haykin2008communication}).}
 \begin{equation}
\label{transmitionRate}
q\triangleq\textbf{Pr}\left\{\gamma_k=1 | |x_k|>d\right\} = 
\textbf{Pr}\left\{\gamma_k=1 | \gamma_k \ne 2 \right\} 
\end{equation} 

\begin{asp}[Processor availability]
\label{processor}
The processor is triggered by arrival of valid data. The processor availability for control at different time-instants is (conditionally) i.i.d. Thus, we denote by  $N_k \in \left\{ 0,1,2,\cdots,N_{max} \right\}$, how many processing units are available at time instants $k$. The process ${\{ N \}}_{\mathbb{N}_0}$ has conditional probability distribution:
\begin{align}
\label{Nkprobability}
&\textbf{Pr}\left\{N_k=j | \gamma_k=1 \right\}=p_j, \ \ j \in \left\{ 0,1,2,\cdots,N_{\max} \right\},
\end{align}
where $p_j \in [0,1)$ are given.	\\
For other values of $\gamma_k$, no plant input is calculated. Thus the processing resources are considered not available regardless, i.e.:
$$
\textbf{Pr}\left\{N_k=0 | \gamma_k \in \{0,2\} \right\}=1
$$
\end{asp}

The next assumption is a combination of some Assumptions from \cite{quevedo2013sequence} and \cite{huang2014event}.
 \begin{asp} [Coarse and fine control policy]~
 \label{VkSAC}

The coarse control law $\kappa_1:\mathbb{R}^n \rightarrow \mathbb{U}$ requires $1$ processing unit to compute, whereas the fine control law $\kappa_2:\mathbb{R}^n \rightarrow \mathbb{U}$ requires $\eta \in \mathbb{N},~\eta \geq 1,$ processing units to compute. We also assume that there exist a common Lyapunov function $V : \mathbb{R}^n \rightarrow \mathbb{R}_{\geq 0} $; $\varphi_1,\varphi_2 \in \mathcal{K}_\infty$, and $\rho_1\in \mathbb{R}_{\geq 0}$, $\alpha>0 $, such that $\forall x \in \mathbb{R}^n$  \begin{align}
 &\varphi_1\left(\left|x\right|\right) \leq V(x)\leq \varphi_2\left(\left|x\right|\right) \label{vAssump}
   \\
 &V\left(f\left(x,{0}\right)\right) \leq \alpha V(x) ~~\text{(open-loop bound)}\label{openLoopBound}\\
 &V \left( f \left( x, \kappa_1(x) \right) \right) \leq {\rho_1}V(x)  ~~(\text{closed-loop contraction 1}) \label{rho1}\\
  &V \left( f \left( x, \kappa_2(x) \right) \right) \leq {\rho_2}V(x)  ~~(\text{closed-loop contraction 2}) \label{rho2}
 \end{align}
and the fine control policy $\kappa_2$ is better than the coarse policy $\kappa_1$ in the sense that
$
\rho_2 < \rho_1$. \hfill $\blacksquare$
\end{asp} 
\subsection{Event-triggered sequence-based anytime algorithms with one control policy}
The baseline algorithm, here denoted by $B1$, amounts to a direct implementation of $\kappa_1$ as per 
\begin{equation}
u_k=\left\{\begin{array}{ll}
\kappa_1(x_k)&\text{if}~\gamma_k=1~\text{and processor is available} \\
0 & \text{otherwise}
\end{array} \right. \notag
\end{equation}

\begin{figure}[t]
\begin{center}
\scalebox{0.99}{
\tikzstyle{splitRec}=[rectangle split,rectangle split parts=2,draw,text centered]

\begin{tikzpicture}[auto, node distance=2cm,>=latex']
    % We start by placing the blocks
    \node [coordinate](startPoint) {} ;
    \node [splitRec,rounded corners, below of=startPoint, node distance=1.5cm, align=center] (caseWithinBound) {$\gamma_k=2$ \nodepart[align=left]{second} $b_k\leftarrow \textbf{0}_\Lambda$ \\ OUTPUT $u_k=0$};
    \node [splitRec,rounded corners,above right of =caseWithinBound, node distance=7.5cm, align=center] (caseOutside) {$\gamma_k=1$~\\$(\left|x_k>d\right|\&~success ~transmission)$ \nodepart[align=left]{second} If $N_k=0$ then $b_k\leftarrow \Xi b_{k-1}$\\If $N_k>0$ \\~~~$\chi \leftarrow x_k,~j\leftarrow 1$ \\~~~REPEAT\\~~~~~~Compute $b_k(j)=\kappa_1(\chi)$\\~~~~~~Update $\chi \leftarrow f(\chi,b_k(j)),~j\leftarrow j+1$\\~~~UNTIL $j>N_k$\\ OUTPUT $u_k=b_k(1)$ };
     \node [splitRec,rounded corners,below right of=caseOutside, node distance=7.5cm, align=center] (caseDropout) {$\gamma_k=0$  \nodepart[align=left]{second} $b_k\leftarrow \Xi b_{k-1}$ \\ OUTPUT $u_k=b_k(1)$};
 %\node [coordinate,north west of=caseWithinBound, node distance=1cm](fooLeft) {} ;
    \draw[<->,bend left,line width=1] (caseWithinBound.north west) to (caseOutside.west);
    \draw[<->,bend left,line width=1] (caseOutside.east) to (caseDropout.north east) ;
    \draw[<->,line width=1] (caseWithinBound) to (caseDropout); 
\end{tikzpicture}}
\end{center}   
\caption{Operation of $A1$, $b_k$: buffer content at time instant $k$. }
\label{algoA1}
\end{figure}
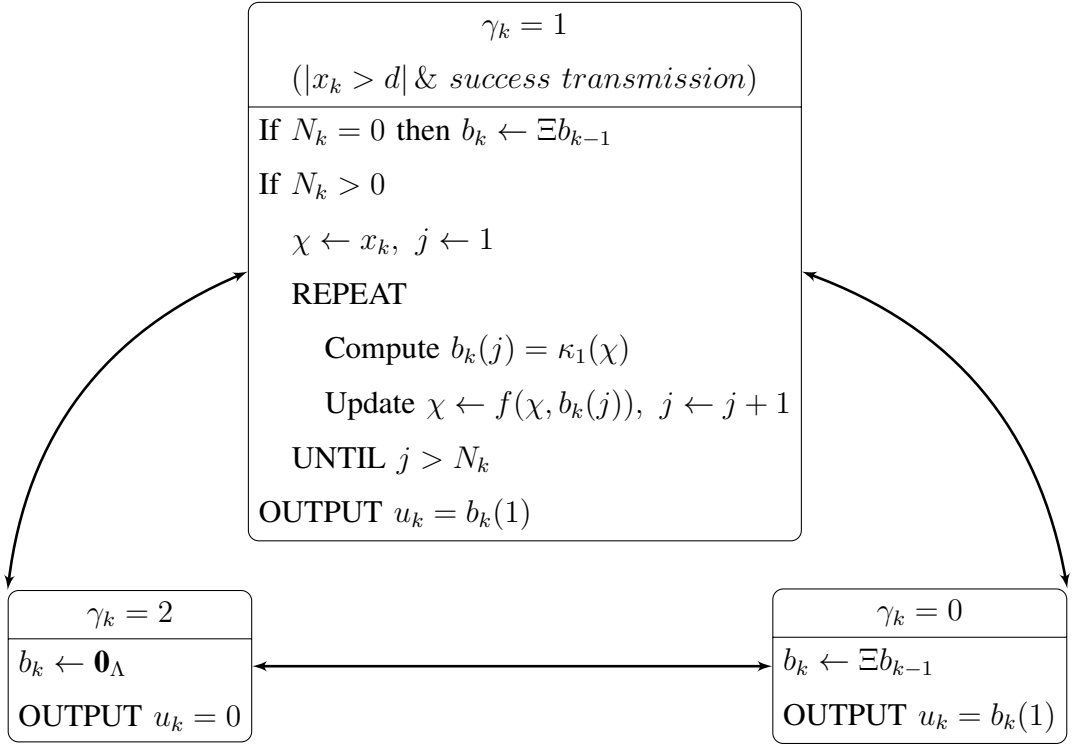

Fig. \ref{algoA1} shows the operation of the (one-control law) E-SAC in \cite{quevedo2014stochastic}. We denote this algorithm as $A1$. In $A1$, tentative future inputs using $\kappa_1$ are calculated and stored in a local buffer $b_k\in \mathbb{R}^{\Lambda}$ ($\Lambda$: buffer size, the maximum number of control inputs it can store), whenever the computing resources are available ($\gamma_k=1~ \&~ N_k >0$). When $|x_k|>d$ and processing resources are unavailable (due to dropouts or unavailable processor, i.e. $\gamma_k=0$ or $N_k=0$), the buffer is shifted, i.e., the first element is thrown away and the rest is kept. If $|x_k|\leq d$ ($\Leftrightarrow\gamma_k=2$), the buffer is cleared. The matrix representing the shift operation is defined as
  \begin{equation}
 \Xi=\left[\begin{matrix}
 {0} &~1~&0&~\cdots ~&0 \\
 0 &~0~&1&~\cdots& ~0 \\
\vdots & \vdots & \vdots & \ddots & \vdots \\
 0 & 0 & 0 & \cdots & 1 \\
 0 & 0 & 0 & \cdots & 0\\
 \end{matrix} \right] \in \mathbb{R}^{\Lambda  \times \Lambda } \notag
 \end{equation}

The first element in the buffer is used for the current input.
We also refer $A1$ to as {\it one-control law} E-SAC scheme. For more details, see \cite{quevedo2014stochastic}.

\subsection{E-SAC with multiple control laws}
In control system design, at times one may encounter situations where one would like to switch between different control laws in respond to changing operating conditions. For example, depending on available computational resources, one may switch between a suboptimal or optimal controller, a short or long prediction horizon MPC, or a fixed- or floating-point controller implementation. In our conference contribution \cite{huang2014event}, E-SAC was extended to schemes featuring two control laws, $\kappa_1$ and $\kappa_2$, to capture such situations. 
We refer to $\kappa_1$ as the coarse (baseline) control law and to $\kappa_2$ as the fine control law. The fine control law $\kappa_2$ requires more computational resource to execute than $\kappa_1$ as shown in Assumption \ref{VkSAC}.    

\subsubsection{Algorithm \textbf{$B2$}: two-control law E-SAC without buffer}
Algorithm \textbf{$B2$} amounts to a direct implementation of $\kappa_1$ and $\kappa_2$ without any buffering. The plant input is calculated as
\begin{equation}
\notag
\ u_k= \left \{
\begin{array}{lll}
\kappa_2 \left(x_k\right), \ \ &\text{if}~ N_k \geq \eta  			 \\
\kappa_1 \left(x_k\right), 	&\text{if}~	0 < N_k < \eta				\\
{0}, 	& \text{if}~ N_k=0.
\end{array} \right.
\end{equation}

\subsubsection{Algorithm \textbf{$A2$} two-control law E-SAC with buffer}
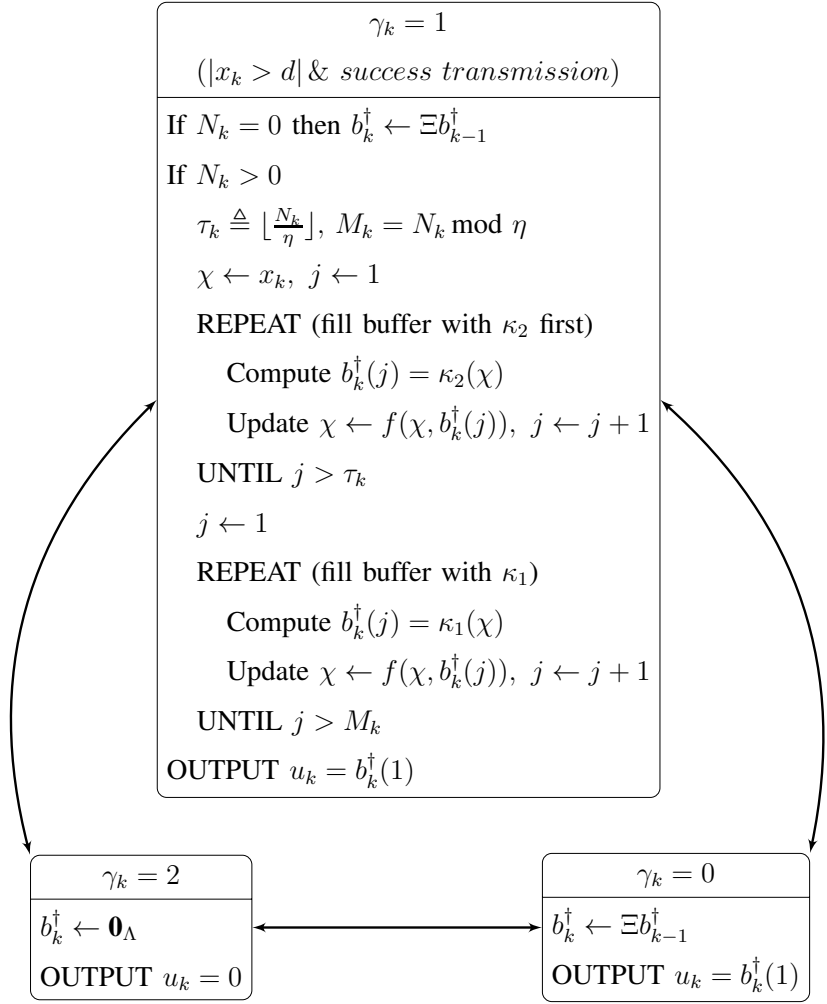
\begin{figure}[t]
\begin{center}
\scalebox{0.9}{
\tikzstyle{splitRec}=[rectangle split,rectangle split parts=2,draw,text centered]

\begin{tikzpicture}[auto, node distance=2cm,>=latex']
    % We start by placing the blocks
\node [coordinate](startPoint) {} ;

\node [splitRec,rounded corners,below of =startPoint, node distance=5.2cm, align=center] (caseOutside) {$\gamma_k=1$\\$(\left|x_k>d\right|\&~success ~transmission)$ 	\nodepart[align=left]{second} If $N_k=0$ then $b^\dagger_k\leftarrow \Xi b^\dagger_{k-1}$\\If $N_k>0$\\~~~$\tau_k\triangleq\lfloor \frac{N_k}{\eta} \rfloor,$ $M_k= N_k\, \mbox{mod } \eta$\\~~~$\chi \leftarrow x_k,~j\leftarrow 1$ \\~~~REPEAT (fill buffer with $\kappa_2$ first)\\~~~~~~Compute $b^\dagger_k(j)=\kappa_2(\chi)$\\~~~~~~Update $\chi \leftarrow f(\chi,b^\dagger_k(j)),~j\leftarrow j+1$\\~~~UNTIL $j>\tau_k$\\~~~$j\leftarrow 1$\\~~~REPEAT (fill buffer with $\kappa_1$)\\~~~~~~Compute $b^\dagger_k(j)=\kappa_1(\chi)$\\~~~~~~Update $\chi \leftarrow f(\chi,b^\dagger_k(j)),~j\leftarrow j+1$\\~~~UNTIL $j>M_k$\\ OUTPUT $u_k=b^\dagger_k(1)$ };
\node [coordinate,right of =caseWithinBound, node distance=3.95cm ](rightOfMidlle) {} ;
\node [coordinate,left of =caseWithinBound, node distance=3.95cm ](leftOfMidlle) {} ;
\node [splitRec,rounded corners, below of=leftOfMidlle, node distance=11.5cm, align=center] (caseWithinBound) {$\gamma_k=2$ \nodepart[align=left]{second} $b^\dagger_k\leftarrow \textbf{0}_\Lambda$ \\ OUTPUT $u_k=0$};
\node [splitRec,rounded corners,below of=rightOfMidlle, node distance=11.5cm, align=center] (caseDropout) {$\gamma_k=0$  \nodepart[align=left]{second} $b^\dagger_k\leftarrow \Xi b^\dagger_{k-1}$ \\ OUTPUT $u_k=b^\dagger_k(1)$};
     
\draw[<->,bend left,line width=1] (caseWithinBound.north west) to (caseOutside.west);
\draw[<->,bend left,line width=1] (caseOutside.east) to (caseDropout.north east) ;
\draw[<->,line width=1] (caseWithinBound) to (caseDropout); 
\end{tikzpicture}}
\end{center}   
\caption{Operation of $A2$, $b^\dagger_k$: buffer content at time instant $k$. }
\label{algoA2}
\end{figure}

Fig. \ref{algoA2} shows the operation of $A2$.
Similar to $A1$, a local buffer with contents $b^\dagger_k$ of size $\Lambda \in \mathbb{N} $ is used to store the sequence of tentative future plant inputs calculated by either $\kappa _1$ or $\kappa _2$ at time $k$ using excess processing resources. 

To be more specific, the control policies $\kappa_1$ or $\kappa_2$ and their tentative future sequence will be executed depending on the values of $N_k$, as illustrated next.

Given the available processing unit $N_k>0$ (assumed known in advance), we can write 
$$ %\begin{equation}
   N_k=\tau_k\eta+M_k
%   \label{tau_M}
$$ %\end{equation}
where $\tau_k\triangleq\lfloor \frac{N_k}{\eta} \rfloor,$ and $M_k= N_k\, \mbox{mod } \eta$ ($\eta$: see Assumption \ref{VkSAC}). Firstly, a tentative control sequence is computed by iterating $\tau_k$ times the model \eqref{sysmodel} using $\kappa_2$, and then by iterating $M_k$ times the model \eqref{sysmodel} using $\kappa_1$. 
 When $|x_k|>d$ and computing resources are unavailable (due to dropouts or unavailable processor), the buffer is shifted. If $|x_k|\leq d~(\Leftrightarrow \gamma_k=2)$, the buffer is cleared. 
The first element in the buffer will be used as the current input. We refer to $A2$ as {\it two-control law} E-SAC scheme. For ease of exposition, in this work, we only present in details the case of two control laws. The case of multiple control laws can be adapted directly.

\begin{remark}
Processing units (in Assumption 1) represent the computational resources (e.g. memory units and given CPU time) available for computing the sequence of predicted control inputs.
 We assume that the processing time of the control task is significantly smaller than the sampling time of the plant model. Here it is important to note that the control values written into the buffer at time $k$ only use information about $x_k$ (if available), but not $x_{k+1} $, or other future states. Merely predictions are used. Therefore, and assuming that processing and transmissions are ``infinitely'' fast, it is appropriate to use a time-invariant system model as \eqref{sysmodel}. The overall system (including communications, and computations) turns out to be stochastically switching, leading to non-trivial dynamics.
\end{remark}

\textbf{Example 1.} Suppose that $N_{\max}=3, \eta=2, \Lambda=3$ and that the processor availability is such that $N_0=3, N_1=0,N_2=2$; the system state is such that $|x_k|>d,~\forall k=0,1,2$ and there are no dropouts.

If algorithm \textbf{$A1$} is used with $\kappa_1$, then the buffer contents become:
$$ \small
\begin{aligned}[c]
&\{b_0,b_1,b_2\}=	
 \left\{
\left[ {\begin{array}{*{20}c}  \kappa_1(x_0) \\ \kappa_1(f(x_0,\kappa_1(x_0))) \\ \kappa_1(\hat{x}_0^{[\kappa_1,2]}) \end{array}} \right]\footnotemark[2],
\left[ {\begin{array}{*{20}c}  \kappa_1(f(x_0,\kappa_1(x_0)))  \\ \kappa_1(\hat{x}_0^{[\kappa_1,2]}) \\ {0} \end{array}} \right], \right.
\left. \left[ {\begin{array}{*{20}c}  \kappa_1(x_2) \\ \kappa_1(f(x_1,\kappa_1(x_2))) \\ {0} \end{array}} \right],
\right\} 
\end{aligned} 
$$
which gives the plant inputs $u_0=\kappa_1(x_0), u_1=\kappa_1(f(x_0,\kappa_1(x_0))), u_2=\kappa_1(x_2)$.
\footnotetext[2]{$\hat{x}_0^{[\kappa_1,2]}=f(f(x_0,\kappa_1(x_0)),\kappa_1(f(x_0,\kappa_1(x_0))))$}

If algorithm \textbf{$A2$} is used, then the buffer contents at times $k \in \{ 0,1,2 \}$ become:
$$ \small
\begin{aligned}[c]
\{b^\dagger_0,b^\dagger_1,b^\dagger_2\}	
= \left\{
\left[ {\begin{array}{*{20}c}  \kappa_2(x_0) \\ \kappa_1(f(x_0,\kappa_2(x_0))) \\ 0 \end{array}} \right],
\left[ {\begin{array}{*{20}c}  \kappa_1(f(x_0,\kappa_2(x_0)))  \\ 0 \\ {0} \end{array}} \right], \right.
\left. \left[ {\begin{array}{*{20}c}  \kappa_2(x_2) \\ 0 \\ {0} \end{array}} \right],
\right\}
\end{aligned}
$$
which gives the plant inputs $u_0=\kappa_2(x_0), u_1=\kappa_1(f(x_0,\kappa_2(x_0))), u_2=\kappa_2(x_2)$.

For the no-buffering schemes, the plant inputs are $u_0=\kappa_1(x_0),~u_1=0,~u_2=\kappa_1(x_2)$ for algorithm $B1$, and $u_0=\kappa_2(x_0),~u_1=0,~u_2=\kappa_2(x_2)$ for algorithm $B2$.

This example suggests that $A2$ outperforms $A1$ since $\kappa_2$ gives better control inputs than $\kappa_1$. The no-buffering schemes $B1$ and $B2$ cannot provide a control input when the processor is unavailable at time step $k=1$.
\subsection{State-dependent random-time drift condition approach for stability analysis of E-SAC}
 In \cite{quevedo2014stochastic}, the state-dependent random-time drift (SRD) condition is developed to derive the stochastic stability of the one-control law scheme with buffering, i.e. the scheme $A1$. For deriving the stability condition, it requires one to calculate probability mass function (pmf) of random variable $\Delta_i$, which denotes the time interval that the buffer becomes empty again. An analytical formulation of this pmf as well as the closed-form for stability boundary of $A1$ has been established in \cite{quevedo2014stochastic}. 
 
 In \cite{huang2014event}, the derivation of stability condition for the two-control law with buffering, i.e. scheme $A2$, follows the same SRD method of \cite{quevedo2014stochastic}. Since there are two control laws in the buffer, the fine control law $\kappa_2$ and the coarse control law $\kappa_1$, there are not only $\Delta_i$ is a random variable, but also the number of times each control law is active, denoted by $r_i$, is also a random variable. Therefore, it is a combinatoric problem and quickly becomes intractable.
 \section{Markov Jump System with Event-triggering Model}
In this section, we propose a different approach for the stability analysis of E-SAC based on Markov jump system ideas. We shall begin our analysis by developing a stochastic model of the buffer contents at any time $k$.
\begin{remark}
\label{a2Reduce2A1}
If $\eta=1$ and $\kappa_2 \equiv \kappa_1$, algorithm $A2$ reduces to $A1$. In addition, $A2$ reduces to $B2$ when the buffer size $\Lambda=1$. Therefore, in the sequel, we only present the stability analysis for algorithm $A2$, since results for $A1$ and $B2$ can be recovered as a special case of $A2$. \hfill $\blacksquare$
\end{remark}
    
 \subsection{Markov state of the buffer content of $A2$}
For two-control law scheme $A2$, we model the content of the buffer via $\theta_k=(F_k;C_k)$ where $F_k$ and $C_k$ indicate the number of $\kappa_2$ (fine control law) and  $\kappa_1$ (baseline or coarse control law) in the buffer at time step $k$ respectively. During periods when $|x_k|>d$, then the transition of $\theta_k$ only depends on $N_k$  which is i.i.d. Hence, during periods when $|x_k|>d$, $\{\theta_k\}$ is a Markov chain. The corresponding state space is 
\begin{align}
\label{stateA2}
&\mathcal{S}^{(A2)}=\{S_i=(\lfloor \frac{i-1}{\eta}\rfloor;(i-1) ~\text{mod} ~\eta)\},
\\
&~i=1,2,...,N_{\max}+1 \notag
\end{align}
and is associated with the conditional ($|x_k|>d$) transition probability  matrix $$
   \Pi^{(A2)}=\{\pi^{(A2)}_{ij}\},~i,j=1,2,...,N_{\max}+1.
   $$ where $\pi^{(A2)}_{ij}={\bf Pr}\{\theta_{k}=S_j|\theta_{k-1}=S_i,\gamma_k\ne 2,\gamma_{k-1}\ne 2\}$ (see Appendix A for the derivation of the entries of the probability transition matrix).
  \begin{figure}[t]
       \begin{center}
       \scalebox{0.85}{
       \begin{tikzpicture}[->, >=stealth', auto, semithick, node distance=3cm]
       \tikzstyle{every state}=[fill=white,draw=black,thick,text=black,scale=1]
       \node[state,line width=3pt,align=center]    (S0)                     {$(0,0)$};
       \node[state,line width=3pt,align=center]    (S1)[right of=S0]   {$(0,1)$};
       \node[state,line width=3pt,align=center]    (S2)[ below right of=S1]   {$(1,0)$};
        \node[state,line width=3pt,align=center]    (S3)[below left of=S2]   {$(1,1)$};
       \node[state,line width=3pt,align=center]    (S4)[below  of=S0]   {$(2,0)$};
       \path
        (S0) edge[loop left]     node[near end]{$l_0$}         (S0)
            edge[bend left]     node[near end]{$l_1$}     (S1)
            edge     node[pos=0.9]{$l_2$}      (S2)
            edge   node[pos=0.8,below]{$l_3$}      (S3)
            edge   node[pos=0.7,left]{$l_4$}      (S4)
        (S1)  edge [in=10,out=35,loop] node{$l_1$}           (S1)
        edge        [above]          node[red]{$l_0$}      (S0)
        edge    [bend left]  node[below,pos=0.75]{$l_2$}      (S2)
                 edge [bend left]  node[near end]{$l_3$}      (S3)
                 edge [bend right]  node[pos=0.8,above]{$l_4$}      (S4)
        
        (S2)  edge  [in=30,out=60,loop]  node[pos=0.95]{$l_2$}           (S2)
             edge  node[red,pos=0.9,above]{$l_0$} (S0)
             edge      node[pos=0.7,above]{$l_1$}                       (S1)
             edge   node[near start,above]{$l_3$}      (S3)
             edge    node[pos=0.9]{$l_4$}      (S4)
     (S3)  edge [loop right] node{$l_3$}           (S3)
               edge[bend left]    node[red,pos=0.2,right]{$l_1+l_0$}                       (S1)
               edge [bend right]  node[above]{$l_2$}      (S2)
               edge  [bend left]node[near end]{$l_4$}      (S4)
      (S4)  edge [loop left] node{$l_4$}           (S4)
                       edge      node[pos=0.9]{$l_1$}                       (S1)
                       edge[bend left]  node[red,pos=0.7,below]{$l_2+l_0$}      (S2)
                       edge  node[below]{$l_3$}      (S3);
      % \node[above=0.5cm] (A){Patch G};
       %\draw[red] ($(D)+(-1.5,0)$) ellipse (2cm and 3.5cm)node[yshift=3cm]{Patch H};
       \end{tikzpicture}}
       \end{center}   
       \caption{State transition diagram for buffer contents of $A2$ ($N_{max}=4,~\eta=2$), $|x_k|>d$. During periods when $|x_k|\leq d$, the buffer is always emptied, i.e. there is no change in the buffer content. }
       \label{markovA2}
         \end{figure}
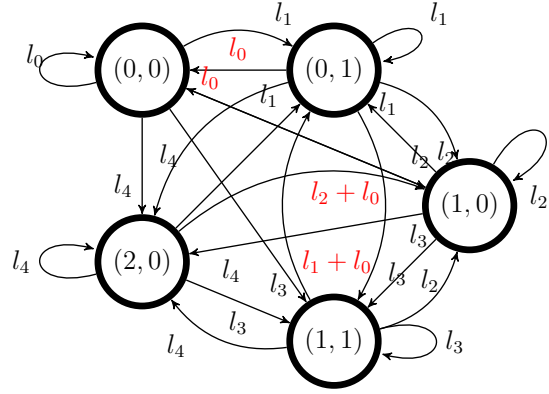
As an example, if $N_{max}=4,~\eta=2$, we have the state space of $A2$: 
 \begin{align*}
 \mathcal{S}^{(A2)}=&\{S_1=(0;0),~S_2=(0;1),~S_3=(1;0),~S_4=(1;1),S_5=(2;0)\}
 \end{align*}
 and the state transition diagram shown in Fig. \ref{markovA2}.

 \begin{remark}
 \label{policyInA2}
In $A2$, for $S_i$ that has the first entry greater than $0$, $\kappa_2$ will be active, i.e., if $\theta_k$ is in the set  $\{S_{\eta+1},S_{\eta+2},...,S_{N_{\max}+1}\}$ we have $F_k>0$, therefore,  $\kappa_2$ is active. If $\theta_k$ is in the set $\{S_2,S_3,...,S_{\eta}\}$, then we have $F_k=0$ and $C_k>0$, therefore, $\kappa_1$ is active. Finally if $\theta_k=S_1=(0;0)$, then a zero control input is used. \hfill $\blacksquare$
 \end{remark}

We note that, when $|x_k|>d$, the buffer and the control input change depending on external factors such as processor and measurement availability which are random. We call this as ``stochastic mode''. On the other hand, when $|x_k|\leq d$, the buffer is cleared and becomes empty. The control, for simplicity and without lost of generality, is set as $u_k=0$. We call this the ``deterministic mode'', since the control is fixed.
\subsection{Markov jump system model}
For our subsequent analysis, it is convenient to introduce
\begin{equation}
\label{ykDef}
y_k\triangleq [x_k;x_{k-1};...;x_{k-\Lambda+1}] \in \mathbb{R}^{n\Lambda}
\end{equation}

\begin{lemma}
\label{aggregatedProcess}
$\{y\}_{\mathbb{N}_0}$ is Markovian.
\end{lemma}
\begin{proof} 
See Appendix B.
\end{proof}

From \eqref{sysmodel} and \eqref{ykDef}, for the process $\{y_k\}$, there exists $\psi: \mathbb{R}^{n\Lambda} \times \mathbb{U} \rightarrow \mathbb{R}^{n\Lambda}$ such that $y_{k+1}=\psi(y_k,u_k)$.

The control $u_k$ is determined by $y_k$ and the random process $\theta_k$ describing the buffer contents (see Remark \ref{policyInA2}). Then, in general, we have $u_k=\hat{\kappa}_{(\theta_k)}(y_k)$
where $\theta_k\in \{S_i\}_{i=1,2,...}$ (see \eqref{stateA2} for $S_i$ ) and set of control laws, associated with the contents of the buffer, $\hat{\kappa}_{(S_i)}: \mathbb{R}^{n\Lambda}\rightarrow \mathbb{U}$. This leads to 
\begin{equation}
\label{ykDynamic}
y_{k+1}=\psi (y_k,\hat{\kappa}_{(\theta_k)}(y_k))
\end{equation}
which is a Markov jump system during intervals when $|x_k|>d$. Here, we use $\hat{\kappa}$ to present the mappings from the domain of $y_k$ (which is different from the domain of $x_k$) to the domain of $u_k$ \footnotemark[3]. \footnotetext[3]{$\kappa_1,\kappa_2$ are the mapping from the domain of $x_k$ to the domain of $u_k$.} 

The mapping from a Markov state to a control law for $\{y_k\}$ process is 
\begin{equation}
\begin{split}
&\hat{\kappa}_{(\theta_k=S_0)}(y_k)=0,\\
&\hat{\kappa}_{(\theta_k=S_1)}(y_k)=\kappa_1(x_k),~\cdots,~\hat{\kappa}_{(\theta_k=S_{\eta)}}(y_k)=\kappa_1(x_k),\\
&\hat{\kappa}_{(\theta_k=S_{\eta+1})}=\kappa_2(x_k),~\cdots,~\hat{\kappa}_{(\theta_k=S_{\Lambda+1})}(y_k)=\kappa_2(x_k)
\end{split}
 \label{mapping}
\end{equation} 
 
\begin{figure}[t]
\scalebox{0.99}{
\begin{tikzpicture}[->, >=stealth', auto, semithick, node distance=3cm]
\tikzstyle{every state}=[fill=white,draw=black,thick,text=black,scale=0.7] 
\node[coordinate] (beginpoint) {};
\node[smallsquare,right of=beginpoint,node distance=1.45cm] (sw1) {};
\node[smallCir,above right of=sw1,node distance=1cm] (sw2) {};
\node[smallCir,below right of=sw1,node distance=1cm] (sw3) {};
\node[ctrlBlk, right of=sw2,node distance=1.9cm] (stochasticCtrl) {$\hat{\kappa}_{(\theta_k)}(y)$};
\node[ctrlBlk, right of=sw3,node distance=1.9cm] (detCtrl) {$\hat{\kappa}_{(-1)}(y)=0y$};
\node[below of= detCtrl,node distance=0.5cm] (detsetup) {\scriptsize \spaceskip0pt \it Deteministic controller (open-loop)};
\node[coordinate, above of=stochasticCtrl,node distance=2.5cm] (markovChain) {$\hat{\kappa}(\theta_k)$};
\node[smallsquare,right of=beginpoint,node distance=6.55cm] (intpoint) {};
\node[ctrlBlk, right of=intpoint,node distance=1.9cm] (plant) {Plant $y_k$};
\node[coordinate,right of=plant,node distance=1.75cm] (endpoint) {};
\node[coordinate,below of=intpoint,node distance=1.7cm] (feedbackpoint) {};
\node[coordinate,right of=sw1,node distance=0.2cm] (foo1) {};
\node[coordinate,above of=foo1,node distance=0.2cm] (foo2) {};
\node[coordinate,below right of=foo1,node distance=0.5cm] (foo3) {};
\node[coordinate,right of=beginpoint,node distance=0.95cm] (foo7) {};
\node[above of=foo7,node distance=0.75cm,align=center] (event1) { \small $\left |y_k(1) \right |> d$\\ {\scriptsize ($\beta_k=1$)}};
\node[below of=foo7,node distance=0.95cm,align=center] (event2) { \small $\left |y_k(1) \right |\leq d$ \\ {\scriptsize ($\beta_k=0$)}};
\node[right of=markovChain,node distance=3.1cm] (legen) {Markov switching};
\draw[dashed] (markovChain) circle (1.55cm);
\node[coordinate,below of=markovChain,node distance=1.55cm] (mcpoint){};
\node[below right of=mcpoint,node distance=0.4cm] (theta){$\theta_k$};
\node[below of= stochasticCtrl,node distance=0.5cm,align=left] (stoSetup) {\scriptsize \it Stochastic controller}; 
%\draw[->] (stochasticCtrl) -- (stoSetup);
\node[state,above of= markovChain,node distance=1.7cm]    (S1)                     {$S_2$};
\node[state,below left of= S1,node distance=1.6cm]    (S0)                     {$S_1$};
\node[state,below right of= S1,node distance=1.8cm]    (S2)                     {$S_3$};
\node[state,below of= S2,node distance=1.3cm]    (S3)                     {$\cdots$};
\node[state,below of= S0,node distance=1.7cm]    (S4)                     {$S_{...}$};

\node[coordinate,right of=beginpoint,node distance=-0.1cm] (boundint){};
\node[coordinate,above of=boundint,node distance=1.3cm] (bound1){};
\node[coordinate,right of=bound1,node distance=2.5cm] (bound2){};
\node[coordinate,below of=bound2,node distance=2.4cm] (bound3){};
\node[coordinate,left of=bound3,node distance=2.5cm] (bound4){};
\draw[dashed,-] (bound1)--(bound2)--(bound3)--(bound4)--(bound1);
\node[coordinate,right of=bound1,node distance=0.9cm] (mid){};

\node[above of=mid,node distance=0.3cm] (legenbound){\scriptsize Event-triggering};

\node[smallCir, right of=stochasticCtrl,node distance=1.8cm] (swS) {};
\node[smallCir, right of=detCtrl,node distance=1.8cm] (swD) {};
\draw [->,line width=2pt](intpoint) to  (swS);
\draw[-] (stochasticCtrl) --(swS);
\draw[-] (detCtrl) --(swD);
\node[coordinate,left of=swS,node distance=0.3cm] (fooS){};
\node[coordinate,above of=fooS,node distance=0.3cm] (boundS1){};
\node[coordinate,right of=boundS1,node distance=1.2cm] (boundS2){};
\node[coordinate,below of=boundS2,node distance=2cm] (boundS3){};
\node[coordinate,left of=boundS3,node distance=1.2cm] (boundS4){};
\draw[-,dashed] (boundS1)--(boundS2)--(boundS3)--(boundS4)--(boundS1);
\node[coordinate,right of=boundS1,node distance=0.6cm] (legendS1){};
\node[above of=legendS1,node distance=0.2cm] (legendS2){};
\node[coordinate,left of= intpoint,node distance=0.4cm] (apivot){};
\node[coordinate,above right of= apivot,node distance=0.2cm] (swsync1){};
\node[coordinate,below of= apivot,node distance=0.4cm] (swsync2){};
\draw [bend right,<->] (swsync1) to (swsync2);

\draw [bend left,->] (S0) to (S1);
\draw [bend left,->] (S1) to (S2);
\draw [->] (S0) to (S2);
\draw [->] (S0) to (S3);
\draw [->] (S1) to (S0);
\draw [bend right,<->] (S1) to (S3);
\draw [<->] (S2) to (S4);
\draw [->] (S2) to (S1);
\draw [bend right,->] (S2) to (S3);
\draw [bend right,->] (S3) to (S2);
\draw [bend right,->] (S3) to (S4);
\draw [bend right,->] (S4) to (S3);
\draw [bend right,->] (S0) to (S4);

%\draw[-] (stochasticCtrl) -|(intpoint);
%\draw[-] (detCtrl) -|(intpoint);
\draw[->] (intpoint) -- node{$u_k$}(plant);
\draw[-] (plant) -- node{$y_k$}(endpoint);
\draw[-] (endpoint) |-(feedbackpoint);
\draw[-] (feedbackpoint) -|(beginpoint);
\draw[->] (beginpoint) --(sw1);
\draw[->,line width=2pt] (sw1) --(sw2);
\draw[->] (sw2) --(stochasticCtrl);
\draw[->] (sw3) --(detCtrl);
\draw[bend left,<->] (foo2) to (foo3);
\draw[->,line width=1.2pt] (mcpoint) to (stochasticCtrl);
\end{tikzpicture}}
\caption{E-MJS model of E-SAC schemes.}
\label{event_trigger_R2}
\end{figure}
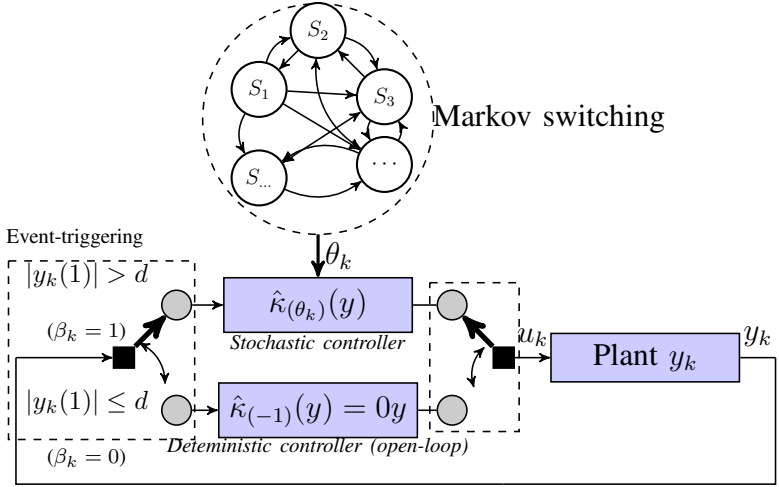

Figure. \ref{event_trigger_R2} shows an equivalent model of the E-SAC schemes $A2$ via the process $\{y_{k}\}$. We call this model the event-triggered Markov jump system model (E-MJS). We use $\beta_k \in \{1,0\}$ to represent the threshold-based triggering event  $|x_k|>d$ ($\gamma_k\in \{0,1\}$) and $|x_k|\leq d$ ($\gamma_k=2$). In this model, the particular schemes such as $A1$ and $A2$ are encoded in the state space $\{S_1,S_2,...\}$ (for e.g. \eqref{stateA2}). The dropouts and processor availability are encoded in the transition probabilities of the state space.
\subsection{General model} 
We now propose a general mathematical description for the E-MJS model of the E-SAC. Consider a non-linear system $y_{k+1}=\psi(y_k,u_k)$ controlled by two controllers: (1) stochastic controller and (2) deterministic controller. The loop is closed with either the stochastic or the deterministic controller (Fig. \ref{event_trigger_R2}).
When the triggering condition is met, the stochastic controller will be deployed. We use $\beta_k\in \{0,1\}$ to indicate the triggering event at time $k$:
\begin{align*}
\beta_k=\left\{ \begin{array}{ll}
0 &\text{use deterministic  controller} \\
1 &\text{use stochastic controller} 
\end{array} \right.
\end{align*}

\subsubsection{Stochastic controller} Due to the external environment, such as time-varying processing powers or dropouts in the communication channels, the controller switches stochastically within a set of $M$ %designed 
 control laws  $\{\hat{\kappa}_{(i)},~ i=1,2,...,M\}$. In this case, the closed loop system model is
\begin{align*}
%\label{sys}
y_{k+1}=\psi(y_k,u_k)=\psi(y_k,\hat{\kappa}_{(\theta_k)}(y_k)) 
\end{align*} where $\{\theta_k\}_{\mathbb{N}}$ is a discrete Markov chain with state space 
\begin{equation}
\label{ssMarkov}
\mathcal{S}=\{1,2,...,M\}
\end{equation}
and the (conditional) transition probability matrix 
$$
\Pi=\{\pi_{ij}\},~~{i,j=1,...,M}
$$
where
$
%\label{piIjDef}
\pi_{ij}={\bf Pr}\{\theta_{k+1}=j|\theta_k=i,\beta_k=1,\beta_{k+1}=1\}
$

Here, for ease of notation, we use numeric representation for the state space of the Markov chain in \eqref{ssMarkov} instead of using $S_i$ as \eqref{stateA2}.

\subsubsection{Deterministic controller}
This controller gives a fixed control policy $u_k=\hat{\kappa}_{(-1)}(y_k)$ and in this case, the closed loop is
$
%\label{sys2}
y_{k+1}=\psi(y_k,\hat{\kappa}_{(-1)}(y_k)).
$

For simplicity but without loss of generality, we set $\kappa_{(-1)}(y_k)=0$.

\section{Stochastic stability of E-MJS model}
\label{section_stability}
In this section, we derive the stochastic stability condition for the proposed E-MJS model.

First, we shall make the following assumptions:
\begin{asp}
\label{lyapunovCandidate}
There exists a non-negative function $\tilde{V}:~\mathbb{R}^{n_y}~\rightarrow \mathbb{R}_{\geq 0}$ ($n_y$ is dimension of $y$) and coefficients  $\sigma_{(i)} \in \mathbb{R}_{\geq 0},~i=-1,1,2,...,M$ such that 
\begin{align}
\label{lyapFunc}
&\tilde{V}(\psi(y,\hat{\kappa}_{(i)}(y)) )\leq \sigma_{(i)} \tilde{V}(y),~~ \forall i=-1,1,2,..,M 
\end{align}
\end{asp}

\begin{asp}
\label{boundVDeter}
There exists a constant $D \in \mathbb{R}_{+}$ such that $\tilde{V}(y_k)\leq D$, if the deterministic controller setup is in operation. \hfill $\blacksquare$
\end{asp}
\begin{remark}
Assumption \ref{lyapunovCandidate} characterises each control law $\hat{\kappa}_{(i)}$ by a scalar $\sigma_{(i)}$, and bounds the rate of increase of $\tilde{V}(y)$ when a control law $\hat{\kappa}_{(i)}$ is active. In Section V.A we show that Assumptions \ref{lyapunovCandidate} and \ref{boundVDeter} are satisfied in the E-SAC schemes $A1$ and $A2$, whenever Assumptions \ref{VkSAC} is satisfied. However, Assumptions \ref{VkSAC} is potentially conservative as a common Lyapunov function is required. \hfill $\blacksquare$
\end{remark}

Let
$
\tilde{V}_{k}=\tilde{V}(y_k)
$, then we obtain the following stochastic model for $\{\tilde{V}_k\}_{\mathbb{N}_0}$:  
\begin{align*}
\left\{ \begin{array}{ll}
\tilde{V}_{k+1}\leq \sigma_{(\theta_k)}\tilde{V}_k,~~(\underbrace{\theta_k\in \{1,2,...,M\}}_{\text{\tiny }Markov~ jump}) &\text{if}~\beta_k=1 \\
\tilde{V}_{k+1}\leq \sigma_{(-1)} \tilde{V}_k &\text{if}~\beta_k=0 
\end{array} \right.
\end{align*}
or in a compact form as: \begin{align} \label{vkModel}
&\tilde{V}_{k+1}\leq \sigma_{(\theta_k)}\tilde{V}_k,~(\theta_k\in \underbrace{\textcolor{black}{\{-1\}}}_{\text {\tiny deterministic ctrl. }}\cup~\underbrace{\textcolor{black}{\{1,2,...,M\}}}_{\text{\tiny Markov~ jump, stochastic ctrl.}})
\end{align}
Note that we have extended the range of $\theta_k$ to include
$
\theta_k=-1$ ($\text{if}~\beta_k=0
$, deterministic mode) to have the compact form as shown in \eqref{vkModel}.

\begin{theorem} 
\label{mstable}
If there exist positive real numbers $\nu_{(1)},...,\nu_{(M)}$,  and $\zeta_{(1)},...,\zeta_{(M)}$ such that 
\begin{align}
\label{mssTheorem1}
\sum\limits_{j=1}^{M}\pi_{ij}\sigma_{(i)}\zeta_{(j)}-\zeta_{(i)}=-\nu_{(i)} ,
\end{align}
for all $i=1,...,M$, then 
$\textbf{E}\{\tilde{V}_k\} <C_1\xi^k\textbf{E}\{\tilde{V}_0\}+C_2 <\infty $
where $$\xi=1-\frac{\min_{1\leq j\leq M}\{\nu_{(j)}\}}{\max_{1\leq k\leq M} \{\zeta_{(k)} \}}\in (0,1)$$ $C_1=\frac{\zeta_{\max}}{\zeta_{\min}},~C_2=\frac{1}{(\zeta_{\min})(1-\xi)}(\max \{ \zeta_{\min}D, |\zeta_{\max}\sigma_{(-1)}-\xi \zeta_{\min}| D \})$ ($\zeta_{\max}=\max\{\zeta_{(i)}\}_{i=1,...,M}$, $\zeta_{\min}=\min\{\zeta_{(i)}\}_{i=1,...,M}$).
\end{theorem}
\begin{proof}
The proof is essentially an adaptation of the general stability result of Markov jump linear systems from \cite{ji1990jump},\cite{fang2002stochastic} 
specialised for \eqref{vkModel}, which is a scalar and positive system with event-triggering. See details in Appendix C.
\end{proof}

Theorem \ref{mstable} provides a general condition for stochastic stability of $\{\tilde{V}_k\}$ in terms of the boundedness property of the expectation. Note that \eqref{mssTheorem1} represents a system of linear equations and can be represented as
\begin{align}
\label{lyaEq}
(I-\Phi\Pi)\underline{\zeta}=\underline{\nu}
\end{align}
where
\begin{align}
\label{energyMatrix}
\Phi\triangleq\text{diag}\{\sigma_{(1)},\sigma_{(2)},...,\sigma_{(M)}\}
\end{align}
and $\underline{\zeta}\triangleq(\zeta_{(1)}~~\zeta_{(2)}~...~\zeta_{(M})^T$, $\underline{\nu}\triangleq(\nu_{(1)}~~\nu_{(2)}~...~\nu_{(M)})^T$. \\
Then, %the condition in 
Theorem \ref{mstable} can be restated as follows:
\begin{corol}%lemma}
\label{alterThereom1}
Define the certification matrix $\mathcal{T}=\Phi \Pi$. If $\mathcal{T}$ is Schur stable, then 
$\textbf{E}\{\tilde{V}_k\} <C_1\xi^k\textbf{E}\{\tilde{V}_0\}+C_2 <\infty $
where $\xi \in (0,1)$, $C_1,C_2 \in \mathbb{R}_{>0}$.
\end{corol} %{lemma}

{\bf Proof:} See Appendix D.

\section{Stochastic Stability of E-SAC Schemes}
In this section, we derive stochastic stability conditions for the E-SAC schemes by applying the results of Section~\ref{section_stability}.

\subsection{Existence and bounds of $\tilde{V}(y_k)$ of E-SAC}
For the process $\{y_k\}$ describing E-SAC (see \eqref{ykDef}), we choose the following function  
\begin{align*}
\tilde{V}:\mathbb{R}^{n\Lambda}\rightarrow \mathbb{R}_{\geq 0},~~\tilde{V}(y_k)= V([I~~0~~0~~...~0]y_k)=V(x_k)
\end{align*}
where $V(.)$ is defined as in Assumption \ref{VkSAC}.

The reason for choosing this $\tilde{V}$ is that it allows us to obtain the bound $\sigma_{(i)}$ in \eqref{lyapFunc}. This bound is related to a control law $\hat{\kappa}_{(i)}$ which is associated with a Markov state. 

Assumptions \ref{VkSAC} and equations \eqref{mapping}, lead to the following bounds for $\tilde{V}(y_k)$ 
\begin{align}
 \label{rhoA2}
&A2:\left\{\begin{array}{l}
\tilde{V}(\psi(y_k,\hat{\kappa}_{(S_1)}(y_k))\leq \alpha\tilde{V}(y_k) \\
\tilde{V}(\psi(y_k,\hat{\kappa}_{(S_i)}(y_k))\leq \rho_1\tilde{V}(y_k),~i=2,3,...,\eta \\
\tilde{V}(\psi(y_k,\hat{\kappa}_{(S_i)}(y_k))\leq \rho_2\tilde{V}(y_k),\\~~~i=\eta+1,...,N_{max}+1 \\
\end{array} \right.
\end{align}
where $\alpha,~\rho_1$ and $\rho_2$ are defined as in Assumptions \ref{VkSAC}.

To show that Assumption \ref{boundVDeter} is also satisfied, we recall that in the deterministic controller setup, $|x_k|\leq d$. Therefore $\tilde{V}(y_k)=V(x_k)\leq \varphi_2(|x_k|)\leq \varphi_2(d)\triangleq D$ (see Assumption \ref{VkSAC}). \hfill $\blacksquare$

\subsection{Stochastic stability for E-SAC schemes $A1$ and $A2$}
We need the following Lemma to establish closed-loop stability, when $A1$ or $A2$ are used.
\begin{lemma}
\label{schurR2}
Consider a $2\times 2$ block matrix $H=\left[ \begin{matrix} X& Y \\Z& M \end{matrix} \right]$, where $X \in \mathbb{R}^{1\times 1}, Y \in \mathbb{R}^{1\times m}, Z \in \mathbb{R}^{m\times 1}, M \in \mathbb{R}^{m\times m} $, and $M$ is Schur stable with non-negative entries and $||M||_{\infty}<1$ and $trace(M^2)<1$. Then $H$ is Schur stable if and only if $g(1)>0$ where $g(\lambda)=(\lambda I_n -X)-Y(\lambda I_m-M)^{-1}Z$.
\end{lemma}

{\bf Proof:} See Appendix E.\hfill $\blacksquare$

Closed-loop stability when using algorithm $A2$ is then established as follows:
\begin{corol}[Stochastic stability of $A2$]
\label{corA2stability}
~

The E-SAC scheme $A2$ yields a stochastically stable loop, in the sense that $V(x_k)$ satisfies the bound condition in Theorem \ref{alterThereom1}, if the certification matrix $${\mathcal{T}}^{(A2)}\triangleq\text{diag}\{\alpha,\underbrace{{\rho_1},...,{\rho_1}}_{\eta-1~\text{times}},\underbrace{{\rho_2},...,{\rho_2}}_{N_{max}-\eta+1~\text{times}}\}\Pi^{(A2)}$$ is Schur stable.

Further, if $\rho_1,\rho_2 <1$ the Schur stability of the certification matrix ${\mathcal{T}}^{(A2)}$ reduces to
 \begin{equation}
 \label{stableA2cond}
 \Psi \triangleq l_0\alpha +l_0\alpha\Theta_2^T(I-G_{{\rho_1},{\rho_2}})^{-1}E_2 <1
 \end{equation}
where $E_2^T=[\rho_1\underbrace{0 ~\cdots~ 0}_{\eta-2~\text{elements}}~\rho_2~\underbrace{0~\cdots~0}_{N_{max}-\eta+1~\text{elements}}]$, $\Theta_2^T=[l_1~~l_2~ \cdots ~l_{N_{max}}]$, and\\
    $G_{{\rho_1},{\rho_2}}={\mathcal{T}}^{(A2)}_{2:(N_{max}+1),2:(N_{max}+1)}$ is the lower right block of ${\mathcal{T}}^{(A2)}$.
\end{corol}
{\bf Proof:} Appendix F. \hfill $\blacksquare$
\begin{remark}
Corollary \ref{corA2stability} provides an analytical expression for the stability boundary of $A2$ which has not been obtained in the earlier works  \cite{huang2014event} or \cite{quevedo2014stochastic}. It is reassuring that  \eqref{stableA2cond} agrees with numerical results of \cite{huang2014event}.
\end{remark}
\begin{remark}
We see that using our approach, the condition that both $\rho_1$ and $\rho_2$ are strictly less than $1$ is not necessary (which was needed in \cite{huang2014event}).  
\end{remark}
\subsection{Recover stability of $A1$}
As aforementioned in Remark \ref{a2Reduce2A1}, $A2$ reduces to $A1$ when $\eta=1$ and $\kappa_2 \equiv \kappa_1$. The probability transition of buffer content of $A2$ when $\eta=1$ (i.e. $A1$ this case) is shown in Appendix as the matrix $\Pi^{(A1)}$ in \eqref{tranMatrixR2}. From Corollary \ref{corA2stability}, we obtain the stability for $A1$:
\begin{corol}[Stochastic stability of $A1$]
\label{corR2Stability}
~

The E-SAC scheme $A1$ yields a stochastically stable loop, in the sense that $V(x_k)$ satisfies the bound condition in Theorem \ref{alterThereom1}, if ${\mathcal{T}}^{(A1)}\triangleq \text{diag}\{\alpha,\underbrace{{\rho_1},...,{\rho_1}}_{N_{max}-1~\text{times}}\} \Pi^{(A1)}$ is Schur stable.

Futher, if $\rho_1<1$ the Schur stability condition for ${\mathcal{T}}^{(A1)}$ is \begin{equation}
\label{stabilityR2}
\Omega\triangleq l_0\alpha(1+ \rho_1 \Theta^T(I-\rho_1 G)^{-1}E_1)<1
\end{equation}
where $\Pi^{(A1)}$ is as \eqref{tranMatrixR2},
$\Theta^T=[l_1~ l_2~ \cdots~l_{N_{max}}]$, $E_1=[1~0~0~\cdots~0]^T\subset \mathbb{R}^{N_{max} \times 1}$, $G=\Pi^{(A1)}_{[2:(N_{max}+1);2:(N_{max}+1)]} \subset \mathbb{R}^{N_{max} \times N_{max}} $ is the lower right block of $\Pi^{(A1)}$ obtained by eliminating the first row and the first column.
\end{corol}

\begin{remark}
Interestingly, \eqref{stabilityR2} is the same stability condition as already derived with different method in \cite{quevedo2014stochastic}. 
\end{remark}
\begin{remark}
The stability of $A1$ is independent of the triggering threshold $d$ as showed in \cite{quevedo2014stochastic}. Similarly, the stability of $A2$ showed in Corollary \ref{corA2stability} is also independent of $d$. The threshold $d$ does however determine the size of the region that the system state converges to. In detail, in Theorem 1, the size of this region is $C_2=\frac{1}{(\zeta_{\min})(1-\xi)}(\max \{ \zeta_{\min}D, |\zeta_{\max}\sigma_{(-1)}-\xi \zeta_{\min}| D \})$. For E-SAC schemes $A1$ and $A2$, $D=\varphi_2(d)$ (as defined in Section V.A.) influences $C_2$.
\end{remark}

\begin{figure}[t!]
\centering
\hspace*{-0.05cm}\includegraphics[scale=0.455]{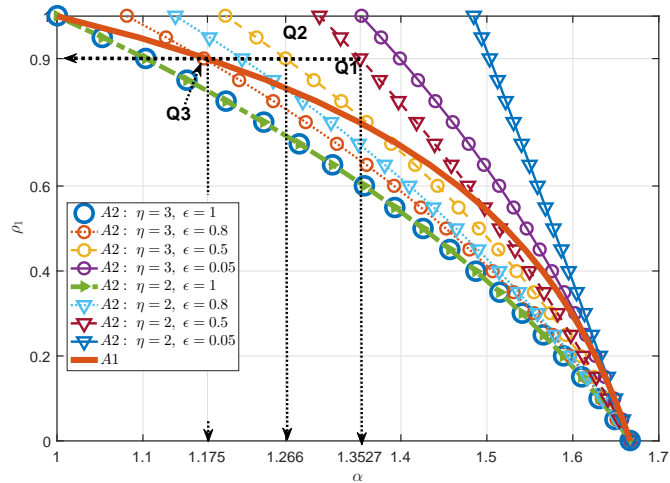}
\caption{ Stability guarantee region for $A2$ with respect to $\eta$ (processing units needed to compute $\kappa_2$) and $\epsilon=\frac{\rho_2}{\rho_1}$. Two-control scheme $A2$ with configuration $Q1$ guarantees to yield a stable system if the open-loop bound satisfying $\alpha<1.3527$. Two-control scheme $A2$ with configuration $Q2$ guarantees to yield a stable system if $\alpha<1.266$. One-control scheme $A1$ with configurationn $Q3$ guarantees to yield a stable system if $\alpha<1.175$. System \eqref{example} has open-loop bound $\alpha=1.35$.  } 
\label{stabilityA2eta23}
\end{figure}
\begin{figure}[t!]
\centering
\hspace*{-0.25cm}\includegraphics[scale=0.469]{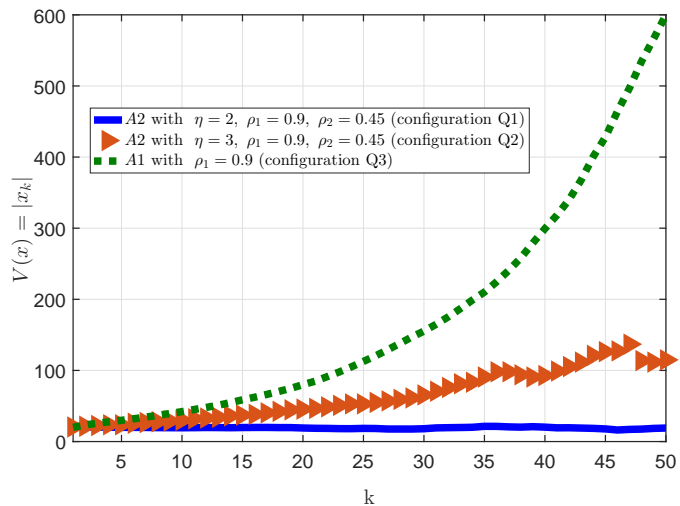}
\caption{Averaged value of $V(x)=|x_k|$ over $10^4$ random realizations.}
\label{trajectoryXk}
\vspace*{-0.25cm}
\end{figure}
\section{Numerical Simulation}
We assume a plant with dynamics
\begin{equation}
x_{k+1}=-1.34x_k+0.01sin(x_k)+u_k+w_k,~~x_0=20 
\label{example}
\end{equation}
where the disturbance $w_k$ is i.i.d., normally distributed with zero mean and unit variance. 
For the proposed schemes with two control laws, we adopt
\begin{align}
&\kappa_1(x_k)=1.34x_k-0.01sin(x_k)+0.9|x_k|,	\label{exK1}	\\
&\kappa_2(x_k)=1.34x_k-0.01sin(x_k)+c_2|x_k|.\label{exK2}
\end{align} 
where $c_2$ is decided later. 

By choosing $V(x)=|x|$, we obtain the open-loop bound $\alpha=1.35$, closed-loop contractions $\rho_1=0.9$ in \eqref{rho1} and $\rho_2=c_2$ in \eqref{rho2}.

We also assume that the buffer size $\Lambda=4$ and that the maximum available processing units are $N_{max}=4$.

The probability of successful transmission is given by:
$$
  q = \textbf{Pr}\{\gamma_k=1 | \ |x_k|>d\}= 0.5
$$

For $N_k \in \{ 0,1,2,\cdots,4 \}$, we assume that the probabilities $p_j=\textbf{Pr}\{N_k=j | \gamma_k=1\}$ are equal for each  $j \in \{ 0,1,2,\cdots,4 \}$, i.e., $p_0=p_1=...=p_4=0.2$.

For ease of presenting the stability region, we introduce a parameter $\epsilon$ that represents the ratio between the closed-loop contractions of control laws $\kappa_1$ and $\kappa_2$ :
$$
\epsilon=\frac{\rho_2}{\rho_1}
$$
We see that $\epsilon \in [0,1]$. It can be said that the smaller the $\epsilon$ is, the ``better'' the second control law $\kappa_2$ is. 

Fig. \ref{stabilityA2eta23} shows the region for the open-loop bound ($\alpha$) and closed-loop contractions ($\rho_1~ \&~\rho_2=\epsilon\rho_1$) that $A2$ and $A1$ guarantee to yield a stochastically stable system, i.e. the region is represented by eq. \eqref{stableA2cond} and \eqref{stabilityR2} given that $\rho_1\in [0,1]$. For a specific value of $\rho_1$ and $\epsilon$, one can figure out the maximum open-loop bound $\alpha$ that is allowed so that the system is guaranteed to be stochastically stable.  The region when using $A2$ depends on parameters $\eta$ and $\epsilon$. Point $Q1=(1.3525,0.9)$ is on the curve with label ``$A2:~\eta=2,~\epsilon=0.5$'' (dashed line with triangle). This implies that given two control laws: (1) the coarse control law as eq. \eqref{exK1} ($\rho_1=0.9$) and (2) the fine control law as eq. \eqref{exK2} with parameters $\eta=2$ and $c_2=\rho_2=0.5*0.9=0.45$, the algorithm $A2$ yields a stochastically stable system {\bf if} the open-loop bound $\alpha$ satisfies $\alpha < 1.3527$. Since system \eqref{example} has $\alpha=1.35<1.3527$, $A2$ yields a stochastically stable system with this configuration of two control laws. The blue line in Fig. \ref{trajectoryXk} confirms this, as the averaged value of $V(x_k)$ is bounded. 

Similarly, by looking at point $Q2$, it shows that $A2$ with $\rho_1=0.9,~\rho_2=0.5*0.9=0.45,~\eta=3$ only guarantees to yield a stochastically stable system if the open-loop bound $\alpha$ satisfies $\alpha<1.266$. And by looking at point $Q3$ on the curve labelled as $A1$, it shows that the one control law $A1$ with $\rho_1=0.9$ (i.e. using control law as \eqref{exK1}) only guarantees to be able to stochastically stabilise a system if the open-loop bound $\alpha$ satisfies $\alpha<1.175$.  Indeed, the averaged value of $V(x_k)$ can be very large in these two cases, see Fig. \ref{trajectoryXk} the triangle points and dotted line, since the open-loop bound of \eqref{example} is $\alpha=1.35$ bigger than the allowed open-loop bounds of these two configurations.
\begin{remark}
Eq. \eqref{stableA2cond} and \eqref{stabilityR2} are sufficient conditions for stochastic stability of $A2$ and $A1$, respectively. Currently, necessary conditions are not available for these schemes.
\end{remark}

\section{Conclusion}
We propose a general model and a novel stability analysis method for event-triggered sequence-based anytime schemes based on Markov jump systems ideas. The proposed method is, unlike the State-dependent Random-time Drift condition approach, scalable for more sophisticated schemes. It also allows us to obtain an analytical expression for the stability boundary of two-control law schemes, as well as recover the existing stability results of one-control law. Future work and extensions are Markovian processor/sensor availability scenarios, processor scheduling, and the appearance of process noise and model uncertainty.

\section*{Appendix}
\subsection{Probability transition in $A2$}
Firstly, we define
\begin{flalign}
&l_j \triangleq \textbf{Pr}\{ N_k=j | \gamma_k \ne 2,\gamma_{k-1} \ne 2 \}=\textbf{Pr}\{ N_k=j | \gamma_k \ne 2 \}&\notag
  \\
       &=\textbf{Pr}\{N_k=j|~\gamma_k =1\} \textbf{Pr}\{\gamma_k =1|\gamma_k \ne 2\}   +\textbf{Pr}\{N_k=j|~\gamma_k =0\} \textbf{Pr}\{\gamma_k =0|\gamma_k \ne 2\} & \label{ljDef}
     \end{flalign}
satisfy $l_0  = p_0q+(1-q),~ l_j   =p_jq ~~(j=1,2,...,N_{\max})\notag 
$, where $p_j$ and $q$ are defined as in \eqref{transmitionRate} and \eqref{Nkprobability}.

We obtain
 \begin{align*}
&\pi_{ij}={\bf Pr}\{\theta_{k}=S_j|\theta_{k-1}=S_i,\gamma_k \ne 2,\gamma_{k-1} \ne 2\} \\
&\pi_{i1}=l_0,~i\in \{1,2,\eta+1 \},~\pi_{i1}=0~\forall i\ne \{1,2,\eta+1\}\\
&\pi_{ij}=l_{j-1}+l_0,~\forall~j=i-\eta,~i=\eta+2,...,N_{max}\\
&\pi_{ij}=l_j~ \text{otherwise}
  %&\pi{ij}=l_j,~\forall~j=i-1,~i \text{mod} \eta >0,~i=1,2,...,\Lambda\\
  \end{align*}
Then, the probability transition matrix is as follows:
 \begin{align*}
  \Pi^{(A2)}=&
 \left[ \begin{matrix}
 l_0&l_1&l_2&\cdots&\cdots&l_{N_{max}}\\
 l_0&l_1&l_2&\cdots&\cdots&l_{N_{max}}\\
 l_0&l_1&l_2&\cdots&\cdots&l_{N_{max}}\\
 0&l_1+l_0&l_2&\cdots&\cdots&l_{N_{max}}\\
 \ldots&\ldots&\ddots&\ldots&\ldots&l_{N_{max}}\\
 0&l_1&\ldots&l_{N_{max}-2}+l_0&\cdots&l_{N_{max}}
 \end{matrix} \right],
 \end{align*}
when $\eta =2$ as an example.

When $\eta=1$ and $\kappa_2 \equiv \kappa_1$, $A2$ reduces to $A1$, and the corresponding probability transition matrix is
 \begin{align}
 \label{tranMatrixR2}
 \Pi^{(A1)}=\left[\begin{matrix}
  l_0&l_1&l_2&l_3&\cdots&l_{N_{\max}}\\
  l_0&l_1&l_2&l_3&\cdots&l_{N_{\max}}\\
  0&l_1+l_0&l_2&l_3&\cdots&l_{N_{\max}}\\
  0&l_1&l_2+l_0&l_3&\cdots&l_{N_{\max}}\\
  \vdots&\vdots&\vdots&\ddots&\ddots&\vdots\\
   0&l_1&\cdots&\cdots&l_{N_{\max}-1}+l_0&l_{N_{\max}}
  \end{matrix} \right]
 \end{align}   
 
 \subsection{Proof of Lemma \ref{aggregatedProcess}}
 If $N_k>0$, then $u_k$ is determined by the current state $x_k$. If the processor is not available, then $u_k$ has
 been determined by the states which are at most $\Lambda$ time stages old, or is given by $u_k=0$. Since the processor availability is independent of the state, the stochastic process $\{y\}_{\mathbb{N}_0}$ is Markovian.
 \subsection{Proof of Theorem \ref{mstable}}
 We define $\zeta_{\max}=\max\{\zeta_{(i)}\}_{i=1,...,M}$, $\zeta_{\min}=\min\{\zeta_{(i)}\}_{i=1,...,M}$, i.e., we have
 \begin{equation}
 \label{boundQtheta}
 0<\zeta_{\min}\leq \zeta_{(\theta_k)}\leq \zeta_{\max},~~\forall ~\theta_k=1,...,M.
 \end{equation}
 Consider $J_k\triangleq J(y_k,\theta_k)=\zeta_{(\theta_{k})}\tilde{V}(y_k)=\zeta_{(\theta_{k})}\tilde{V}_k$
 where $\theta_k$ is as in Eq. \eqref{vkModel}.
 
We will next establish the drift condition  $$
  \textbf{E}\{ J_{k+1}|y_k,\theta_{k}\} \leq \bar{D}+ \xi J(y_k,\theta_k)$$
where $\bar{D}$ is a constant and is derived later (after eq. \eqref{Ziterative}) in the following.
   
Since $J(y_k,\theta_k)$ becomes undefined when $\theta_k=-1$ (deterministic mode, $\zeta_{(\theta_k=-1)}$ does not exist), without loss of generality, we take (extended value) $\zeta_{(-1)}=\zeta_{\min}$.
 
 If ${\theta_k=-1}$, then $\beta_k=0$, i.e., the deterministic mode is active. We have $\tilde{V}(y_k)\leq D$ (Assumption \ref{boundVDeter}), and we also have $\zeta_{(\theta_{k+1})}\tilde{V}_{k+1}\leq \zeta_{\max}\sigma_{(-1)}\tilde{V}_k$. From $\tilde{V}(y_k)\leq D$, we obtain $\zeta_{\max}\sigma_{(-1)}\tilde{V}_k \leq \xi \zeta_{\min} \tilde{V}_k+|\zeta_{max} \sigma_{(-1)} - \xi \zeta_{\min}| D \leq \xi J(y_k,\theta_k)+|\zeta_{max} \sigma_{(-1)} - \xi \zeta_{\min}| D $. Then it follows that 
 \begin{align}
 &\textbf{E}\{ J(y_{k+1},\theta_{k+1})|y_k,\theta_{k}=0\}
 =\textbf{E}\{ \zeta_{(\theta_{k+1})}\tilde{V}_{k+1}|y_k,\theta_{k}=0\} \notag\\
& \leq \xi J(y_k,\theta_k)+|\zeta_{max} \sigma_{(-1)} -\xi \zeta_{\min}| D  \label{boundTheta0} 
 \end{align}
 
 If ${\theta_k=i>0}$, then $\beta_k=1$. By denoting $\Upsilon_k\triangleq(y_k,\theta_k=i>0)$, using the law of total expectation we have:
 \begin{flalign}
 &\textbf{E}\{ J(y_{k+1},\theta_{k+1})|y_k,\theta_{k}=i>0\}=& \notag 
 \\
 &\textbf{E}\{ J(y_{k+1},\theta_{k+1})|\Upsilon_k,\beta_{k+1}=0\}\textbf{Pr}\{\beta_{k+1}=0|\Upsilon_k\}+& \notag
 \\
 &
 \textbf{E}\{ J(y_{k+1},\theta_{k+1})|\Upsilon_k,\beta_{k+1}=1\}\textbf{Pr}\{\beta_{k+1}=1|\Upsilon_k\}& \notag\\
 &\leq \textbf{E}\{ J(y_{k+1},\theta_{k+1})|\Upsilon_k,\beta_{k+1}=0\} +\textbf{E}\{ J(y_{k+1},\theta_{k+1})|\Upsilon_k,\beta_{k+1}=1\}& \label{boundThetak11}
 \end{flalign}
 For $\beta_{k+1}=0$, this implies that $\theta_{k+1}=-1$. And by Assumption \ref{boundVDeter}, we have $\tilde{V}_{k+1}\leq D$, then
  \begin{align}
 &\textbf{E}\{ J(y_{k+1},\theta_{k+1})|\Upsilon_k,\beta_{k+1}=0\}=\textbf{E}\{ \zeta_{(-1)}\tilde{V}_{k+1}|\Upsilon_k,\beta_{k+1}=0\}\leq  \zeta_{\min} D,~~~~(\zeta_{(-1)}=\zeta_{\min})\label{boundThetak12}
 \end{align}
 
 For the case $\beta_{k+1}=1$, the stochastic controller is deployed, i.e., $\theta_{k+1}>0$, we also have $\beta_k=1$, since $\theta_k>0$.\\ 
 We have, using the law of total expectation and using \eqref{mssTheorem1},
 \begin{flalign}
 &\textbf{E}\{ J(y_{k+1},\theta_{k+1})|\Upsilon_k,\beta_{k+1}=1\} &\notag\\
 &=\sum\limits_{j=1}^{M}\textbf{Pr}\{\theta_{k+1}=j|\theta_k=i,y_k,\beta_k=1,\beta_{k+1}=1\}.& \notag\\
  &~~~~\textbf{E}\{\tilde{V}_{k+1}\zeta_{(\theta_{k+1})}|\theta_{k+1}=j,y_k,\theta_{k}=i,\beta_k=1,\beta_{k+1}=1\} &\notag\\
 &\leq\sum\limits_{j=1}^{M}(\pi_{i j})(\sigma_{(i)}\tilde{V}_k\zeta_{(j)})=\tilde{V}_k\left(\zeta_{(i)}-\nu_{(i)}\right) ~~\text{(using \eqref{mssTheorem1})}\notag \\
 &=\tilde{V}_k\zeta_{(i)}\left(1-\frac{\nu_{(i)}}{\zeta_{(i)}} \right)=J(y_k,\theta_k) \left(1-\frac{\nu_{(i)}}{\zeta_{(i)}} \right)\leq  \xi J(y_k,\theta_k)& \label{boundThetak13}    %\left(1-\frac{\min_{1\leq j\leq M}\{\nu_{(j)}\}}{\max_{1\leq k\leq M} \{Q(k) \}}\right) \notag \\
 %& \leq  \xi J(y_k,\theta_k) \label{boundThetak13}
 \end{flalign}
 since $$\left(1-\frac{\nu_{(i)}}{\zeta_{(i)}} \right) \leq 1-\frac{\min_{1\leq j\leq M}\{\nu_{(j)}\}}{\max_{1\leq k\leq M} \{\zeta_{(k)} \}}\triangleq \xi, ~\forall i=1,..,M $$
 and $0<\xi<1 $ due to $\zeta_{(k)}>\nu_{(k)}>0,~\forall k=1,2,...,M$ (from \eqref{mssTheorem1}).\\
 Expressions \eqref{boundThetak11}-\eqref{boundThetak13} lead to:
 \begin{align}
 \textbf{E}\{ J(y_{k+1},\theta_{k+1})|y_k,\theta_{k}>0\}\leq \zeta_{\min} D + \xi J(y_k,\theta_k) \label{boundThetak11All}
 \end{align}
 
 From \eqref{boundTheta0} and \eqref{boundThetak11All}, we obtain
 \begin{align*}
&\textbf{E}\{ J(y_{k+1},\theta_{k+1})|y_k,\theta_{k}\} \leq \bar{D}+ \xi J(y_k,\theta_k) 
 \end{align*}
 i.e., we have 
 %\begin{align}
 %&\textbf{E}\{ J(y_{k+1},\theta_{k+1})|x_k,\theta_{k}\} \leq \bar{D}+ \xi J(y_k,\theta_k) 
 %\end{align}
 \begin{align}
 &\textbf{E}\{ J(z_{k+1})|z_k\} \leq \bar{D}+ \xi J(z_k) 
 \label{Ziterative}
 \end{align}
 where
 \begin{align*}
 &z_k\triangleq (y_k,\theta_k)\\
 &\bar{D}=\max \{ \zeta_{\min}D, |\zeta_{\max}\sigma_{(-1)}-\xi\zeta_{\min}| D \}
 \end{align*}
 From \eqref{Ziterative} and using the Markovian property of $\{z_k\}$, \footnote{Since $\{y_k\}$ and $\{\theta_k\}$ are Markovian, see Lemma \ref{aggregatedProcess} and \eqref{vkModel}.} we have 
 \begin{flalign}
 &\textbf{E}\{ J(z_{1})|z_0\} \leq \bar{D}+ \xi J(z_0) & \label{z0in}\\
 &\textbf{E}\{ J(z_{2})|z_1\}=\textbf{E}\{ J(z_{2})|z_1,z_0\} \leq \bar{D}+ \xi J(z_1) \label{z1in}
 \end{flalign}
 %From \eqref{z1in} we have
 %\begin{flalign}
 %&\textbf{E}\{ J(z_{2})|z_1,z_0\} \leq \bar{D}+ \xi J(z_1) \label{z1z0in}
 %\end{flalign}
 Taking expectation $\textbf{E}\{. |z_0\}$ of both sides of \eqref{z1in}, using tower property of expectation and using \eqref{z0in} we have 
 \begin{flalign}
 &\textbf{E}\{ J(z_{2})|z_0 \}=\textbf{E}\{\textbf{E}\{ J(z_{2})|z_1,z_0\}|z_0 \}\leq \bar{D}+ \xi \textbf{E}\{J(z_1)|z_0\} \leq \bar{D}+\xi \bar{D}+\xi^2 J(z_0) & \notag
 \end{flalign}
 By iterating the above procedure, we obtain
 $$
 \textbf{E}\{ J(z_k)|z_0\}\leq \bar{D} (\sum\limits_{i=0}^{i=k-1} \xi^i) +\xi^k  J(z_0)\leq \frac{\bar{D}}{1-\xi} +\xi^k  J(z_0)
 $$
 i.e., $$\textbf{E}\{ J(y_k,\theta_k)|y_0,\theta_0\}\leq \frac{\bar{D}}{1-\xi} +\xi^k  J(y_0,\theta_0)
 $$
 Using the law of total expectation \footnote{Here we assumed  $\theta_0$ is known.}, if $y_0$ is a discrete random variable we have
 \begin{align}
 &\textbf{E}\{ J(y_k,\theta_{k})\}=\sum\limits_{\forall (y_0,\theta_0)}^{}\textbf{E}\{ J(y_k,\theta_{k})|y_0,\theta_0\}\textbf{Pr}\{y_0,\theta_0\} \notag \\
 &\leq \sum\limits_{\forall (y_0,\theta_0)}^{} \left(\frac{1}{1-\xi}\bar{D} +\xi^k  J(y_0,\theta_0) \right) \textbf{Pr}\{x_0,\theta_0\} \notag \\
 &~~=\frac{1}{1-\xi}\bar{D}+ \sum\limits_{\forall (y_0,\theta_0)}^{} \left( \xi^kJ(y_0,\theta_0) \right)\textbf{Pr}\{y_0,\theta_0\} \notag\\
 &~~=\frac{1}{1-\xi}\bar{D}+\xi^k \textbf{E}\{J(y_0,\theta_0)\} \label{lawofTotal}
 \end{align}
  If $y_0$ is a continuous random variable, we have\footnote{\textbf{pdf}: probability density function.} 
  \begin{align}
   &\textbf{E}\{ J(y_k,\theta_{k})\}=\int\limits_{y_0}^{}\textbf{E}\{ J(y_k,\theta_{k})|y_0,\theta_0\}\textbf{pdf}\{y_0,\theta_0\}dy_0 \notag \\
   &\leq \int\limits_{y_0}^{} \left(\frac{1}{1-\xi}\bar{D} +\xi^k  J(y_0,\theta_0) \right) \textbf{pdf}\{y_0,\theta_0\}dy_0 \notag\\
   &~=\frac{1}{1-\xi}\bar{D}+ \int\limits_{y_0} \left( \xi^kJ(y_0,\theta_0) \right)\textbf{pdf}\{y_0,\theta_0\} dy_0  \notag \\
   &~=\frac{1}{1-\xi}\bar{D}+\xi^k \textbf{E}\{J(y_0,\theta_0)\} 
   \label{lawofTotal2}
 \end{align}
   
 Since $\zeta_{max}\tilde{V}_k\geq J(y_k,\theta_{k})=\zeta_{(\theta_k)}\tilde{V}_k\geq \zeta_{min}\tilde{V}_k$, then $\textbf{E}\{ J(y_k,\theta_{k})\} \geq \textbf{E}\{\tilde{V}_k\} \zeta_{\min} $ and $\textbf{E}\{J(y_0,\theta_0)\}\leq \zeta_{\max}\textbf{E}\{\tilde{V}_0\}$. Then, from \eqref{lawofTotal} (or \eqref{lawofTotal2}) we have
 $$
 \zeta_{\min}\textbf{E}\{\tilde{V}_k\}\leq \zeta_{\max}\xi^k\textbf{E}\{\tilde{V}_0\}+\frac{1}{1-\xi}\bar{D}
 $$

\subsection{Proof of Corollary 1}
 Since $\mathcal{T}$ is Schur stable, we have $
 (I-\mathcal{T})^{-1}=I+\sum\limits_{i=1}^{\infty}(\mathcal{T})^i=I+\bar{\mathcal{T}}
 $
 where $\bar{\mathcal{T}}=\sum\limits_{i=1}^{\infty}(\mathcal{T})^i$. As all entries of $\mathcal{T}$ are non-negative, all entries of $\bar{\mathcal{T}}$ are non-negative. Thus, for any given $\underline{\nu}\succ 0$, the solution of \eqref{lyaEq} is $\underline{\zeta}=(I-\mathcal{T})^{-1}\underline{\nu}=\underline{\nu}+\bar{\mathcal{T}}\underline{\nu}$ satisfying $\underline{\zeta} \succ 0$. Then, by applying Theorem 1, we obtain the desired bound for $\textbf{E}\{\tilde{V}_k\}$. \hfill $\blacksquare$
    
 \subsection{Prove of Lemma \ref{schurR2}  }
 Consider the characteristic polynomial $P(\lambda)\triangleq\text{det}(\lambda I-H)=\left| \begin{matrix} \lambda  - X& -Y \\ -Z& \lambda I_m-M \end{matrix} \right|$. We denote $eig(\cal{M})$ as an eigenvalue of a matrix $\cal{M}$. 
 
 For $\lambda$ not inside the unit circle $|\lambda|\geq 1$, we have $eig(\lambda I_m - M)=\lambda-eig(M) \ne 0$ (due to Schur stability of $M$) hence $\lambda I_m - M$ is invertible. Using the determinant result for a 2-by-2 block  matrix, we have 
 \begin{align}
 P(\lambda)&=\text{det}(\lambda I_m-M)\text{det}\left(\lambda - X-Y(\lambda I_m-M)^{-1}Z\right) \notag\\
 &=h(\lambda)g(\lambda),~ ~\forall |\lambda| \geq 1 \label{gLambda}
 \end{align} 
 where $h(\lambda)=\text{det}(\lambda I_m-M)$ and $g(\lambda)=\text{det}\left(\lambda - X-Y(\lambda I_m-M)^{-1}Z\right)$.
 
 ``$\Rightarrow$'' We have $H$ is Schur stable, now we need to proof $g(1)>0$. 
 
 Since $H$ is Schur stable, from the Schur-Cohn criteria (see \cite{lasalle1986stability}, p.27) we have $P(1)>0$. It follows that $h(1)g(1)>0$. As $M$ is Schur stable, we have $h(1)=det(I_m-M) >0$. Hence, $g(1)>0$.  
 
 ``$\Leftarrow$'' We have $g(1)>0$, now we need to show that $H$ is Schur stable. Assume that $H$ is not Schur stable. Then there exists $\lambda_0$, $|\lambda_0| \geq 1$ such that $P(\lambda_0)=0.$ We then have $\lambda_0I_m-M$ is invertible and then $P(\lambda_0)=h(\lambda_0)g(\lambda_0)=0$ where $h(\lambda_0)=det(\lambda_0I_m-M) \ne 0$. It follows that $g(\lambda_0)=0$. 
 
 (1) If $\lambda_0$ is a real number: we can proof that $g(\lambda)$ is a strictly increasing function on $\mathbb{R}$. Then we have $0=g(\lambda_0)>g(1)>0$. This is a contradiction.
 
 (2) If $\lambda_0$ is a complex number, we have $\lambda_0^*$, the conjugate of $\lambda_0$, is also an eigenvalue of $H$. Then $\lambda_0^2$ and $(\lambda_0^*)^2$ are eigenvalues of $H^2$. Then, $trace(H^2)=\sum (eig(H^2))\geq \lambda_0^2+(\lambda_0^*)^2 \geq 2 $. It can also be shown that $trace(H^2)=X^2+2YZ+trace(M^2)<(X+YZ)^2+trace(M^2) <2$, since $g(1)>0$ $\Rightarrow $ $X+YZ<1$ and $trace(M^2)<1$ as in the assumption of Lemma 2. This is also a contradiction.
 
 The two above contradictions establish the result. \hfill $\blacksquare$
     
\subsection{Proof of Corollary \ref{corA2stability}}
  \begin{proof}
The certification matrix of $A2$, $\mathcal{T}^{(A2)}$, is calculated as follows.
From \eqref{rhoA2} and \eqref{energyMatrix}, we have $\Phi^{(A2)}=\text{diag}\{\alpha,\underbrace{{\rho_1},...,{\rho_1}}_{\eta-1~\text{times}},\underbrace{{\rho_2},...,{\rho_2}}_{N_{max}-\eta+1~\text{times}}\}$. 
We then have ${\mathcal{T}}^{(A2)}\triangleq\Phi^{(A2)}\Pi^{(A2)}$. Then, the first statement of Corollary \ref{corA2stability} is directly proved by Corollary \ref{alterThereom1}.

To prove the second statement of Corollary \ref{corA2stability}, firstly we have ${\mathcal{T}}^{(A2)}$ can be written as:\\ 
${\mathcal{T}}^{(A2)}=\left[\begin{matrix}
 l_0\alpha& \alpha \Theta^T\\
 l_0\rho_1 E_1&\rho_1 G \end{matrix} \right]$.

It is easy to prove that the (2,2) block of ${\mathcal{T}}^{(A2)}$, which is $G_{{\rho_1},{\rho_2}}$, is Schur stable by proving that $||G_{{\rho_1},{\rho_2}}||_{\infty}<1$, hence $1> ||G_{{\rho_1},{\rho_2}}||_{\infty}\geq \sigma(G_{{\rho_1},{\rho_2}})$ (spectral radius is less  than or equal any matrix norm). It can also be shown that $trace(G_{{\rho_1},{\rho_2}}^2)<1$. Therefore, by applying Lemma \ref{schurR2} we have the Schur condition for ${\mathcal{T}}^{(A2)}$ is that
$$1-l_0\alpha -l_0\alpha\Theta_2^T(I-G_{{\rho_1},{\rho_2}})^{-1}E_2 >0$$
which gives \eqref{stableA2cond}.
\end{proof}
\balance  
\bibliographystyle{IEEEtran}
\bibliography{refTAC171999}
 %Bayes rule:
 %\begin{align*}
 %&{\bf Pr} \{\Delta_i=j|\beta_{k_i+1}\ne 2\}{\bf Pr} \{\beta_{k_i+1}\ne 2\}\\
 %&={\bf Pr} \{\beta_{k_i+1}\ne 2|\Delta_i=j\}{\bf Pr} \{\Delta_i=j\} \leq \frac{{\bf E}\{V(x_{k_i})\}}{D}\\
 %&{\bf Pr} \{\}
 %\end{align*}
\end{document}